\documentclass{lmcs}
\pdfoutput=1
\usepackage[utf8]{inputenc}

\usepackage{lastpage}
\lmcsdoi{21}{2}{1}
\lmcsheading{}{\pageref{LastPage}}{}{}%
{May~30,~2023}{Apr.~04,~2025}{}

\usepackage{hyperref}
\usepackage{mathtools}
\usepackage{booktabs}
\usepackage{hyperref}
\usepackage{mathtools}
\usepackage{amssymb}
\usepackage{amsmath}
\usepackage{amsthm}
\usepackage{stmaryrd}
\usepackage{ebproof}

\newcommand{\Nat}{\mathbb{N}}
\newcommand{\Rea}{\mathbb{R}}
\newcommand{\Bol}{\mathbb{B}}
\newcommand{\Lis}{\mathbb{L}}

\newcommand{\T}[1]{\mathbf{T}_#1}
\newcommand{\I}{\mathbf{I}}
\newcommand{\J}{\mathbf{J}}
\newcommand{\Iclauses}[1]{\big(\forall_{\vec{x}_i}(\vec{A}_i(#1)\to #1\vec{t}_i)\big)_{i<k}}
\newcommand{\Jclauses}[1]{\big(\forall_{\vec{y}_j}(\vec{B}_j(#1)\to #1\vec{s}_j)\big)_{j<m}}

\newcommand{\INclause}[1]{\Big(\forall_{\vec{y}}\big(\bigwedge_{i<k}(\forall_{\vec{x}_i}(\vec{y}\equiv\vec{t}_i\to B_{i}))\to #1\vec{y}\big)\Big)}
\newcommand{\JNclause}[1]{\Big(\forall_{\vec{y}}\big(#1\vec{y}\to\bigwedge_{i<k}(\forall_{\vec{x}_i}(\vec{y}\equiv\vec{t}_i\to B_{i}))\big)\Big)}

\newcommand{\eqdef}{\coloneqq}
\newcommand{\N}[1]{#1^\mathbf{N}}
\newcommand{\NN}[1]{#1^\mathbf{NN}}

\newcommand{\BISH}{\mathrm{BISH}}

\newcommand{\INT}{\mathrm{INT}}

\newcommand{\TOT}{\Leftrightarrow}
\newcommand{\To}{\Rightarrow}

\newcommand{\CM}{\mathrm{CM}}

\newcommand{\Pred}{\mathrm{\mathcal{P}}}

\renewcommand{\Form}{\mathcal{F}}
\newcommand{\TCF}{\mathrm{TCF}}
\newcommand{\CLSN}{\mathrm{CLSN}}
\newcommand{\No}{\textnormal{\textbf{N}}}
\newcommand{\mueq}{\overset{\mu}{\eqdef}}
\newcommand{\nueq}{\overset{\nu}{\eqdef}}
\newcommand{\tyto}{\rightarrowtriangle}

\title[Strong negation in TCF]{Strong negation in the \texorpdfstring{\\}{} Theory of Computable Functionals TCF}
\author[N.~K{\"o}pp]{Nils K{\"o}pp\lmcsorcid{0000-0002-8280-4744}}
\address{Ludwig-Maximilians Universit{\"a}t, Theresienstr. 39, 80333, M{\"u}nchen}
\email{koepp@math.lmu.de}
\author[I.~Petrakis]{Iosif Petrakis\lmcsorcid{0000-0002-4121-7455}}
\address{Università di Verona, Strada le Grazie 15, 37134, Verona}
\email{iosif.petrakis@univr.it}

\begin{document}
\begin{abstract}
\noindent We incorporate strong negation in the theory of computable functionals $\TCF$,
a common extension of Plotkin's PCF and G\"{o}del's system \textbf{T}, by defining
simultaneously strong negation $A^{\No}$ of a formula $A$ and strong negation $P^{\No}$ of a predicate $P$ in $\TCF$.
As a special case of the latter, we get strong negation of an inductive and a coinductive predicate of $\TCF$.
We prove appropriate versions of the Ex falso quodlibet and of double negation elimination for strong negation
in $\TCF$. We introduce the so-called tight formulas of TCF i.e.,
formulas implied by the weak negation of their strong negation, and the relative tight formulas.
We present various case-studies and examples, which reveal the naturality of our definition of strong negation in $\TCF$
and justify the use of TCF as a formal system for a large part of Bishop-style constructive
mathematics.
\end{abstract}
\maketitle
\section{Introduction}
\label{sec:intro}
In constructive (and classical) logic negation $\neg A$ of a formula $A$ is defined in a negative way as
the implication  $A \to \bot$, where $\bot$ denotes ``falsum''. Constructively though, $\neg A$ has a weak behaviour,
as the negation of a conjunction
$\neg(A \wedge B)$ does not generally imply the disjunction $\neg A \vee \neg B$, and similarly the negation of
a universal
formula $\neg\forall_x A(x)$ does not generally imply the existential formula $\exists_x \neg A(x)$.
Even if $A$ and  $B$ are stable i.e., $\neg \neg A \to A$ and $\neg \neg B \to B$, we only get
$\neg(A \wedge B) \vdash \neg \neg (\neg A \vee \neg B)$ and $\neg\forall_x A(x)
\vdash \neg \neg \exists_x \neg A(x)$, where $\vdash$ is the derivability relation in minimal logic.
For this reason, $\neg A$ is called the \textit{weak} negation of $A$, which, according to Rasiowa~\cite{Ra74}, p.~276,
is not constructive, exactly due to its aforementioned behaviour.

In contrast to weak negation, and in analogy to the use of a strong ``or'' and ``exists'' in constructive logic,
Nelson~\cite{Ne49} introduced the \textit{strong} negation
$\sim \hspace{-1.5mm} A$ of $A$.
Following Kleene's recursive realisability, Nelson developed constructive arithmetic with strong negation and showed that
it has the same expressive power as Heyting arithmetic. Within Nelson's realisability,
the equivalences $\mathbin{\sim} (A \wedge B) \TOT \ \mathbin{\sim}A$ $\vee \mathbin{\sim} B$ and $\mathbin{\sim}\forall_x A(x) \TOT \exists_x \mathbin{\sim} A(x)$ are realisable. In axiomatic presentations of constructive logic with
strong negation $(\CLSN)$ (see~\cite{Ra74, Gu77}) these equivalences are axiom schemes. In most formalisations of $\CLSN$,
but not in all (see system \textbf{N4} in~\cite{AN84}), strong negation $\mathbin{\sim} A$ implies weak negation
$\neg A$. Markov~\cite{Ma50} expressed weak negation
through strong negation and implication as follows:
$\neg A := A \to \mathbin{\sim} A.$
Rasiowa~\cite{Ra74} went even further, introducing a \textit{strong} implication
$A \To B$, defined through the standard ``weak'' implication $\to$ and strong negation as follows:
$A \To B := (A \to B) \wedge (\mathbin{\sim} B \to \mathbin{\sim} A).$
The various aspects of the model theory of $\CLSN$ are developed in many papers
(see e.g.,~\cite{Gu77, Ak88, SV08a, SV08b}).

The critique of weak negation in intuitionistic logic and later constructive mathematics $(\CM)$ goes back to Kolmogorov~\cite{kolmo1925}, it motivated Johanssen to develop minimal logic~\cite{johannson1937} and Griss to consider negationless mathematics~\cite{Gr46} altogether.
Despite the use of weak negation in Brouwer's intuitionistic mathematics $(\INT)$ and Bishop's constructive
mathematics $(\BISH)$, both, Brouwer and Bishop, developed a positive and strong approach
to many classically negatively defined concepts. An apartness relation $\neq_X$ on a set $(X, =_X)$ is strong
counterpart to the weak negation $\neg(x =_X x{'})$ of its equality, strong complement $A^{\neq} := \{x \in X \mid
\forall_{a \in A}(x \neq_X a)\}$ of a subset $A$ of $X$ is the positive version of its standard,
weak complement\footnote{In this case $A$ is considered to be an extensional subset of $X$ i.e., it is defined by separation on $X$ with respect to an extensional property on $X$.} $A^c := \{x \in X \mid \forall_{a \in A}(\neg(x =_X a))\}$,
inhabited sets are strong, non-empty sets, and so on. Recently, Shulman~\cite{Sh21} showed
that most of these strong concepts of $\BISH$ ``arise automatically from an ``antithesis'' translation
of affine logic into intuitionistic logic via a ``Chu/Dialectica construction''.
Motivated by this work of Shulman, and in relation to a reconstruction of the theory of sets underlying
$\BISH$ (see~\cite{Pe20, Pe21, Pe22a} and~\cite{PW22,MWP23}),
the second author develops (in a work in progress) a treatment of strong negation within
$\BISH$, arriving in a heuristic method for the definition of these strong concepts
of $\BISH$ similar to Shulman's. For that, the equivalences occurring in the axiom schemata related to strong negation in
$\CLSN$ become the definitional clauses of the recursive definition of $\mathbin{\sim} A$, with respect to the
inductive definition of formulas $A$ in $\BISH$.

This idea of a recursive definition of strong negation is applied here to the formal theory of computable functionals
$\TCF$, developed mainly by Schwichtenberg in~\cite{SW12}. $\TCF$ is a
common extension of Plotkin's PCF and G\"{o}del's system \textbf{T}, it uses, in contrast to Scott's LCF,
non-flat free algebras as semantical domains for the base types, and for higher types its intended model consists
of the computable functions on partial continuous objects (see~\cite{SW12}, Chapter 7).
Furthermore, it accommodates inductively and coinductively defined predicates and its underlying logic is minimal.
$\TCF$ is the theoretical base of the proof assistant Minlog~\cite{TMS} with which the extraction of computational
content of formalised proofs is made possible. A significant feature of $\TCF$ is its relation to $\CM$, as
many parts of $\BISH$ have already been formalised in $\TCF$ and implemented
in Minlog (see e.g.,~\cite{Sc09, Sc22, Sc23}).

A novelty of the definition of strong negation $A^{\No}$ of a formula $A$ of $\TCF$ is its extension, by the
first author, to inductive
and coinductive predicates. Moreover, the analysis of strong negation in the informal system of non-inductive\footnote{Bishop's informal system of constructive mathematics BISH can be seen as a system that does not accommodate inductive definitions. Although the inductive definitions of Borel sets and of the least function space generated by a set of real-valued functions on a set are found in~\cite{Bi67}, Bishop's original measure theory was replaced later by the Bishop-Cheng measure theory that avoids the inductive definition of Borel sets, and the theory of function spaces was never elaborated by Bishop (it was developed by the second author much later). It is natural to consider the extension BISH$^*$ of BISH as Bishop's informal system of constructive mathematics with inductive definitions given by rules with countably many premises.}, constructive
mathematics $\BISH$ given by the second author is reflected and extended here in the study of strong negation within
the formal system $\TCF$, as results and concepts of the former are
translated and extended in the latter.
We structure this paper as follows:
\begin{itemize}
	\item In Section ~\ref{sec: TCF} we briefly describe $\TCF$: its predicates and formulas, its derivations and axioms, Leibniz equality, and weak negation.
	\item In Section ~\ref{sec: SNTCF} strong negation $A^{\No}$ of a formula $A$ of $\TCF$ is defined by recursion on
	the definition of formulas of $\TCF$. This definition is extended to the definition of
	strong negation $P^{\No}$ of a predicate $P$ of $\TCF$ (Definition \ref{defi:StrongNegation}). As a special case of the latter, we get  strong
	negation of an inductive and a coinductive predicate of $\TCF$ .
	\item In Section ~\ref{sec: propSNTCF} we prove some fundamental properties of $A^{\No}$ in $\TCF$. Among others, we show
	an appropriate version of the Ex-falso-quodlibet (Lemma \ref{lem:StrongNegationLem}) and of double negation elimination shift for strong negation
	(Theorem \ref{thm:DoubleStrongNegElim}). We also explain why there is no uniform way to define  $A^{\No}$ (Proposition \ref{prop:nouniform}), and we incorporate
	Rasiowa's strong implication in $\TCF$. Finally, we show that weak and strong negation are classically equivalent (using Lemma~\ref{lem:StrongNegationLem} and Lemma~\ref{lem:classWNtoSN}).
	\item Motivated by the concept of a tight apartness relation in $\CM$, in Section ~\ref{sec: SNTCFtight} we introduce
	the so-called tight formulas of $\TCF$ i.e., the formulas $A$ of $\TCF$ for which the
	implication $\neg(A^{\No}) \to A$ is derivable in minimal logic.
 As in the partial setting of $\TCF$ very few non-trivial formulas are tight, the introduced notion of relative tightness is studied (see Proposition~\ref{prop:reltightbasic} and Lemma~\ref{lem:reltight}).
	\item In Section ~\ref{sec: ex} some important case-studies and examples are considered that reveal the naturality of our main
	Definition \ref{defi:StrongNegation}. Furthermore, they justify the use of $\TCF$ as a formal system for a large part of $\BISH$,
	since the use of strong negation in $\TCF$ reflects accurately the seemingly ad hoc definitions of
	strong concepts in $\BISH$.
\end{itemize}

\section{The Theory of Computable Functionals TCF}
\label{sec: TCF}
In this section we describe the theory of computable functionals $\TCF$.
For all notions and results on $\TCF$ that are used here without further explanation or proof, we refer
to~\cite{SW12}.
\par The term-system of TCF, $\mathbf{T}^+$, is an extension of G\"odels \textbf{T}. Namely, the simple types are enriched with free algebras $\iota$, generated by a list of constructors, as extra base types. Particularly, we have the usual strictly positive data-types, for example
	\[
	\Nat \coloneqq \langle
	\mathbf{0}:\Nat,
	\mathbf{S}:\Nat\tyto\Nat
	\rangle,
	\qquad
	\alpha\times\beta\coloneqq \langle
	[\_\,,\,\_]:\alpha\tyto\beta\tyto\alpha\times\beta\
	\rangle.
\]
Moreover,
nesting into already defined algebras with parameters is allowed e.g., given lists $\Lis(\alpha)$ and boole $\Bol$, the finitely branching boolean-labelled and number-leafed trees is defined via the algebra $\mathbb{T}$ with constructors
\[
    \mathbf{Leaf}:\Nat\tyto\mathbb{T},\quad\mathbf{Branch}:\Lis[\alpha/\mathbb{T}]\times \Bol\tyto\mathbb{T}.
\]
The term-system $\mathbf{T}^+$ is \textit{simple-typed lambda-calculus}. The \textit{constants} are \textit{constructors} or \textit{program constants} $\mathtt{D}:\vec{\rho}\tyto\tau$, defined by a system of consistent left-to-right computation rules. Particularly, we have \textit{recursion} and \textit{corecursion operators} $\mathcal{R}_\iota$ respectively, ${^{co}{\mathcal{R}}}_\iota$, and the destructor $\mathcal{D}_\iota$ for any free algebra $\iota$. In that sense, it is regarded as an extension of \textit{G\"odels} \textbf{T}. In general terms are partial and if needed, computational properties e.g., totality have to be asserted by proof.
\begin{defi}[Types and terms] We assume countably many distinct type-variable-names $\alpha,\beta,\gamma,\xi$ and term-variable-names $x^\rho , y^\rho$ for every type $\rho$. Then we inductively define $\mathcal{Y}$ and $\mathbf{T}^+$ by the following rules:
	\[\begin{aligned}
		\rho,\sigma,\tau\in\mathcal{Y}&\Coloneqq \alpha,\beta,\gamma,\xi\ \big|\ \rho\tyto\sigma\ \big|\ \iota\eqdef\langle\mathtt{C}_i:\vec{\rho}_i[\beta/\iota]\tyto\iota\rangle_{i<k},
		\\[5pt]
		t,s\in\mathbf{T}^+&\Coloneqq\ x^\rho,y^\rho\ \big|\ \left(\lambda_{x^\rho}t^\sigma\right)^{\rho\tyto\sigma}\ \big|\ \left(t^{\rho\tyto\sigma}s^\rho\right)^\sigma\ \big|\ \mathtt{C}^{\vec{\rho_i}\to\iota}_i\ \big|\ \mathtt{D}^{\vec{\rho}\tyto\tau}.
	\end{aligned}
	\]
	\begin{itemize}
		\item For a constructor $\mathtt{C}_i:\vec\rho_{i}[\beta/\iota]\tyto\iota$ we require that $\beta$ occurs at most strictly positive in all its $k<|\vec{\rho}_i|$ argument-types $\rho_{ik}(\beta)$
		\item A program constant $\mathtt{D}:\vec{\rho}\tyto\sigma$ is introduced by $i<k$ computation rules $\mathtt{D}P_i^{\vec{\rho}}:=t_i^{\sigma}$
		, where $P_i$ is a constructor pattern namely, a list of term expressions built up using only variables and constructors.
		\item Consistency is ensured, by requiring that there is no unification for $P_i,P_j$, if $i\neq j$.
		\item $\mathcal{Y}$ and $\mathbf{T}^+$ are closed under substitution and types, respectively terms are viewed as equivalent with respect to bound renaming. Terms are identified according to beta equivalence $(\lambda_xt)s=t[x/s]$ and computation rules.
	\end{itemize}
\end{defi}
\begin{exa}
	We have the following algebras:
	\[
		\begin{aligned}
			\Bol&\coloneqq \langle
				\mathtt{f}:\Bol ,
				\mathtt{t}:\Bol
				\rangle&
			&&\textnormal{(boole)},
			\\
			 \Lis(\alpha)&\coloneqq \langle
				[\,]: \Lis(\alpha),
				\_::\_\,: \alpha\tyto\Lis(\alpha)\tyto\Lis(\alpha)
				\rangle&
			&&\textnormal{(lists)}.
		\end{aligned}
	\]
	For the algebra $\Nat$ we also use roman numerals e.g., $\mathtt{1}=\mathtt{S0}$.
	We introduce a program constant $\mathtt{natbin}:\Nat\tyto\Nat\tyto\Lis[\Bol]$, such that $\mathtt{natbin}\ n\ m$ computes a binary representation of $n+2m$ and it is defined by the following rules:
	\[
			\begin{alignedat}{4}
			\mathtt{nbin}\ \;
			&\mathtt{0}
			&\mathtt{0}
			&\ =\  [\,],
			\\
			\mathtt{nbin}\ \;
			&\mathtt{0}
			&\mathtt{S}k
			&\ =\ \mathtt{f}\, ::\, (\mathtt{nbin}\ \mathtt{S}k\ \mathtt{0}),
			\\
			\mathtt{nbin}\ \;
			&\mathtt{1}
			&k
			&\ =\ \mathtt{t}\, ::\, (\mathtt{nbin}\ k\ \mathit{0}),
			\\
			\mathtt{nbin}\ \;
			&\!\!\mathtt{SS}n\!\!\,
			&k
			&\ =\ \mathtt{nbin}\  n\ \mathtt{S}k.
			\end{alignedat}
	\]
For non-nested algebras with nullary constructor, the \textit{total} terms are finite expressions consisting only of constructors. For the type $\Nat$, total terms are exactly of the form
\[
\mathtt{0},\mathtt{S0},\dots,(\underbrace{\mathtt{S}\cdots\mathtt{S}}_{n\mathrm{\ times}})0,\dots\,.
\]
Similar, \textit{cototal} terms are all finite or infinite constructors-expressions. For the type of unary natural number, the only non-total, cototal term is
\[
    \infty=\mathtt{S}\mathtt{S}\mathtt{S}\cdots,
\]
which has the property $\infty=\mathtt{S}\infty$. For lists of natural numbers, $\mathbb{L}(\Nat)$, there are (countably) many non-total and cototal terms, e.g., constant infinite lists $n:n:\cdots$, or
\[
    1:2:3:4\cdots,\quad 2:3:5:\cdots:(n\mathrm{-th\ prime}):((n+1)(\mathrm{-th\ prime})\cdots.
\]
\end{exa}
To assert properties of programs i.e., \textit{totality} or \textit{correctness}, we need the corresponding formulas. E.g., the intended domain of $\mathtt{nbin}$ above are all finite natural numbers. Formally we introduce a predicate $\T{\Nat}$ with
\[
\mathtt{0}\in \mathbf{T}_\Nat,\qquad \forall_n(n\in \mathbf{T}_\Nat\to ((\mathtt{S}n)\in \mathbf{T}_\Nat)),
\]
and it is the \textit{smallest} predicate with this property i.e.,
\[
	0\in P\to \forall_{n\in\T{\Nat}}(n\in P\to (\mathtt{S}n)\in P) \to \T{\Nat}\subseteq P.
\]
Generally, \textit{formulas} and \textit{predicates} of TCF are structurally parallel to types and algebras, though with additional term-dependency. Namely, we introduce \textit{inductive} and \textit{co-inductive predicates} as \textit{least} respectively, \textit{greatest} fixed-point of clauses i.e., formulas of the form $\forall_{\vec{x}}(A_0\to \dotsb\to A_{k-1}\to X\vec{t})$ which are also called \textit{clauses}.

Formulas are generated by universal quantification $\forall$ and implication $\to$. The logical connectives $\wedge,\vee$ and $\exists$ are special cases of inductive predicates, but here, for the clarity of the presentation, we prefer to introduce them as primitives.
\begin{defi}[Predicates and formulas]
\label{def:predform}
    We assume countably many predicate-variables $X$ for any arity $\vec{\rho}$. We simultaneously define formulas $A,B\in\mathcal{F}$ and predicates $P,Q\in\mathcal{P}$ by the following rules:
	\[
		\begin{aligned}
			\mathcal{P}&\Coloneqq\ X\ \big|\ \{\vec{x}\,|\,A\} \ \big|\ \I\eqdef\mu_X \vec{K} \ \big|\ \J\eqdef\nu_{X}\vec{K},\\
			\mathcal{F}&\Coloneqq\ P\vec{t}\ \big|\ A\diamond B\  \big|_{\diamond\in \{\to,\wedge,\vee\}}\ \nabla_x A\ \big|_{\nabla\in\{\forall,\exists\}}.
		\end{aligned}
	\]
	\begin{itemize}
		\item Sometimes variable-dependencies will be emphasized by writing the variable-name, in parantheses, directly behind an expression e.g., $\nu_Y\vec{K}(Y)$.
		\item $\I\eqdef \mu_X\vec{K}$ denotes the least fix-point closed under the clauses $\vec{K}$.
		\item $\J\eqdef \nu_X\vec{K}$ denotes the greatest fixed-point closed under the clauses $\vec{K}$
        A clause or introductory rule $K_i(X)$ of a fixed-point must be a closed formula of the form $K_i=\forall_{\vec{x}}(\vec{A}_i(X)\to X\vec{t}_i)$, where $X$ occurs at most strictly positive in all premises $A_{ik}(X)$ and $\vec{t}_i:\vec{\rho}$ a list of terms.
		\item We refer to $\{\vec{x}^{\vec{\sigma}}\,|\,A(\vec{x})\}$ , where $\vec{x}$ are exactly the free term-variables of $A$, as a comprehension term. It is a predicate of arity $\vec{\sigma}$ and we identify the formulas $\{\vec{x}\,|\,A\}\vec{y}$ and $A[\vec{x}/\vec{y}]$.
		\item Term and type-substitution extends to predicates. Formulas and predicates are closed under (admissible) simultaneous substitutions for type-, term- and predicate-variables.
        Given two predicates $P,Q$ of the same arity, we use the following abbreviations.
        \[
            \begin{gathered}
                 P\cap Q\eqdef \{\vec{x}\,|\, P\vec{x}\wedge Q\vec{x}\},\qquad
                 P\cup Q\eqdef \{\vec{x}\,|\, P\vec{x}\vee Q\vec{x}\},\\
                 P\subseteq Q\eqdef \forall_{\vec{x}}(P\vec{x}\to Q\vec{x}).
            \end{gathered}
        \]
	\end{itemize}
	Furthermore, for any formula or predicate, we associate the following collections of predicate variables:
	\begin{itemize}
		\item at most strictly positive: $\mathbf{SP}(\cdot)$,
		\item strictly positive and free: $\mathbf{PP}(\cdot)$,
		\item strictly negative and free: $\mathbf{NP}(\cdot)$,
		\item freely occurring: $\mathbf{FP}(\cdot)$.
	\end{itemize}
	These predicate variables are defined recursively by the following rules, where unless one of the specific rules on the right applies, they follow the generic rule for $\mathbf{P}$ on the left:
	\[
	\begin{gathered}
	\begin{aligned}
		&\mathbf{P}(A\diamond B)\eqdef \mathbf{P}(A)\cup\mathbf{P}(B),\\
		&\mathbf{P}(\nabla_x A)\eqdef\mathbf{P}(A),\\
		&\mathbf{P}(\I)\eqdef\sideset{}{_{i<k}}{\bigcup}\mathbf{P}(\vec{A}_i)\backslash \{X\},\\
		&\mathbf{P}(P\vec{t}\;)\eqdef \mathbf{P}(P),\\
		&\mathbf{P}(X)\eqdef\{X\},\\
		&\mathbf{P}(\{\vec{x}\,|\, A\})\eqdef \mathbf{P}(A),\\
            &\mathbf{P}(\vec{A})\eqdef \sideset{}{_{k<|\vec{A}|}}{\bigcup}\mathbf{P}(A_k),
	\end{aligned}\qquad
	\begin{aligned}
		&\mathbf{SP}(A\to B)\eqdef\mathbf{SP}(B)\backslash\mathbf{FP}(A),\\
		&\mathbf{SP}(Y)\eqdef\mathbf{Pvar},\\
		&\mathbf{PP}(A\to B)\eqdef\mathbf{PP}(B),\\
		&\mathbf{NP}(Y)\eqdef\{\},\\
		&\mathbf{NP}(A\to B)\eqdef\mathbf{PP}(A)\cup\mathbf{NP}(B),\\
            &\mathbf{SP}(\I)\eqdef \bigcap\mathbf{SP}(\vec{A}_i)\backslash\{X\},\\
            &\mathbf{SP}(\vec{A})\eqdef \sideset{}{_{k<|\vec{A}|}}{\bigcap}\mathbf{SP}(A_k).
	\end{aligned}
	\end{gathered}
	\]
	Here $\mathbf{Pvar}$ denotes the set of all predicate variables and $\I\eqdef \pi_X\big[\forall_{\vec{x}_i}(\vec{A}_i\to X\vec{t}_i\big]_{i<m}$ ($\pi=\mu,\nu$) is a least- or greatest fixed-point.
\end{defi}
\begin{rem}
	We will refer to predicates introduced as least- and greatest-fixed-point of a list of clauses $\I\eqdef\pi_X\vec{K}$ ($\pi=\mu,\nu$) as \textit{inductive} (respectively \textit{coinductive}) predicate and also write $\I\mueq\vec{K}[X/I]$ (or $\J\nueq\vec{K}(\J)$) to avoid mentioning predicate variables. Assuming this list of clauses if given by
	\[
		\vec{K}=[\forall_{\vec{x}_i}(\vec{A}_i(X)\to X\vec{t}_i)]_{i<k}\subseteq\mathcal{F},
	\]
	then, using the rules and axioms from the next definition, the following equivalence will is derivable
	\[
		\vdash\vec{y}\in\I \leftrightarrow \bigvee_{i<k}\exists_{\vec{x}_i}\big(\vec{y}\equiv\vec{t}_i\bigwedge \vec{A}_i(\I)\big),
	\]
	where $\equiv$ is \textit{Leibniz-equality} which is Definition \ref{def:falsity} below.
        Particularly, we can view the latter expression above as an operator
        \[
            \Phi_\I(X)\eqdef \Big\{\vec{y}\,\Big|\,\bigvee_{i<k}\exists_{\vec{x}_i}\big(\vec{y}\equiv\vec{t}_i\bigwedge \vec{A}_i(X)\big)\Big\}.
        \]
        Since $X$ is at most strictly positive in $\Phi_\I(X)$, it is monotone i.e., $P\subseteq Q\to \Phi_\I(P)\subseteq\Phi_\I(Q)$ (see Proposition \ref{prop:mon}), and $\I$ is a fixed-point of the monotone operator $\Phi_\I$. In practice we introduce inductive and coinductive predicates by supplying a list of clauses, though in light of the above, we will also adapt the very convenient notations:
        \[
            \I\eqdef \mu_X(\Phi(X)\subseteq X),\qquad \J\eqdef\nu_Y(Y\subseteq\Psi(Y)).
        \]
\end{rem}
The logic of TCF is minimal first-order predicate-logic extended with axioms for predicates.
\begin{defi}[Derivations and axioms]
	Derivations are defined according to the rules of natural deduction for minimal logic. For $\wedge$, $\vee$ and $\exists$ we use the following axioms:
	\[
	\begin{aligned}
		&(\wedge)^+:A\to B\to A\wedge B, &&(\vee)^+_0:A\to A\vee B,
		\\
		&(\exists)^+: \forall_x(A\to \exists_x A), &&(\vee)^+_1:B\to A\vee B,
		\\
		&(\wedge)^-:(A\to B\to C)\to A\wedge B\to C,\\
		&(\vee)^-:(A\to C)\to(B\to C) \to A\vee B\to C,\\
		&(\exists)^-:\forall_x(A\to B)\to \exists_x A\to B &&(x\textnormal{\ not\ free\ in\ }B).
	\end{aligned}
	\]
	\noindent Predicates $\I\eqdef\mu_X(\Phi(X)\subseteq X)$ and  $\J\eqdef\nu_Y(Y\subseteq\Psi(Y)$ are endowed with two axioms each, namely:
	\[
		\begin{aligned}
			&(\I)^+: \Phi[X/\I]\subseteq\I,\quad &&(\I)^\mu[P]:(\Phi[\I\cap P]\subseteq P)\to (\I\subseteq P),\\
			&(\J)^-: \J\subseteq\Psi[X/\J], &&(\J)^\nu[Q]:(Q\subseteq(\Psi[\I\cup Q]))\to (Q\subseteq \J).
		\end{aligned}
	\]
	\noindent Assuming  concrete representations
	\[
		\I\eqdef\mu_X\Iclauses{X},\quad\J\eqdef\nu_{Y}\Jclauses{Y},
	\]
	the axioms are of the following form:
	\[
		\begin{aligned}
			&(\I)^+_i:\forall_{\vec{x}_i}(\vec{A}_i(\I)\to \I\vec{t}_i),\qquad
			(\I)^\mu[P]:\big(\forall_{\vec{x}_i}(\vec{A}_i(\I\cap P))\to P\vec{t}\;\big)_{i<k}\to \I\subseteq P,\\[10pt]
			&(\J)^-:\J\vec{y}\to \bigvee_{j<m}\exists_{\vec{y}_j}\big(\vec{y}\equiv\vec{s}_j\wedge\bigwedge \vec{B_{j}}(\J)\big),\\
			&(\J)^\nu[Q]:\forall_{\vec{y}}\Big(Q\vec{y}\to \bigvee_{j<m}\exists_{\vec{y}_j}\big(\vec{y}\equiv\vec{s}_j\wedge\bigwedge \vec{B_{j}}(\J\cup Q)\big)\Big)\to P\subseteq \J.
		\end{aligned}
	\]
	Here $\equiv$ is the Leibniz-equality (defined below) and $\vec{x}\equiv\vec{y}:=\bigwedge_{i}(x_i\equiv y_i)$.
\end{defi}
\begin{rem}
\label{rem:nonrec}
	If $\I\eqdef\mu_X(\Phi(X)\subseteq X)$ is a \textit{non-recursive} inductive predicate i.e., $X$ does not occur freely in $\Phi$, the greatest- and least-fixed-point coincide. Namely, let $\J\eqdef\nu_X(\Phi(X)\subseteq X)$, then the axioms for $\I$ and $\J$ are as follows.
	\[
		\begin{aligned}
			&(\I)^+:\Phi\subseteq\I,\quad &&(\J)^-:\J\subseteq \Phi,\\
			&(\I)^\mu[P]:(\Phi\subseteq P)\to (I\subseteq P),\quad &&(\J)^\nu[Q]:(Q\subseteq\Phi)\to(Q\subseteq \J).
		\end{aligned}
	\]
\end{rem}
\begin{defi}[Leibniz equality]
	\label{def:falsity}
	Leibniz equality $\equiv$ is the inductive predicate
	\[
		\equiv\;\coloneqq \mu_X\forall_x (X\,x x)
	\]
	and codifies the least reflexive relation. Arithmetical falsity and weak negation are directly defined as the following formulae:
	\[
	\bot\eqdef \mathtt{f}\equiv\mathtt{t},\qquad \neg A\eqdef A\to \bot.
	\]
\end{defi}
\begin{rem}
	The two axioms associated to Leibniz equality $\equiv$ are
	\[
	(\equiv)^+:\forall_x x\equiv x,\qquad (\equiv)^\mu[P]:\forall_x Pxx\to \forall_{x,y}(x\equiv y\to Pxy).
	\]
	With these axioms it is immediate to show that $\equiv$ is an equivalence relation. Moreover, we can prove the indiscernibility of identicals
	\[
	\vdash x\equiv y\to A(x)\to A(y).
	\]
	Although we work in minimal logic, we can prove the intuitionistic Ex-falso principle from arithmetical falsity $\bot$. Notice that for any two terms $t,s$ of the same type $\tau$ we can define a program constant $\mathtt{or}_{ts}\;\mathtt{f}=t$, $\mathtt{or}_{ts}\;\mathtt{t}=s$. Assuming $\bot$ then gives $\mathtt{or}_{ts}\;\mathtt{f}\equiv\mathtt{or}_{ts}\;\mathtt{t}$, so
	\[
		\bot\vdash t\equiv s.
	\]
	Then by induction on formula, we can extend \textbf{EfQ} to (almost) all formulas. Namely, we need to exclude some non-terminating inductive predicates e.g., $\mu_X(X\to X)$. For simplicity, we restrict ourselves to inductive predicates $\I$
	that are inhabited i.e., there exists a term $t$ with $t\in\I$. In particular, this is the case if $\I$ has a clause $(i<k)$ and $X$ does not occur freely in any of the premises $A_{i\nu}$.
 \end{rem}

 In the following Proposition we demand all inductive predicates occurring in a formula to be inhabited. As we are using Minimal Logic, Ex-falso with $\bot\coloneqq 0=1$ is a Theorem and it fails for non-inhabited inductive predicates $\I$, since these are literally empty i.e., even assuming $0=1$, it is impossible to prove $t\in\I$, for any term $t$. In practice, non-inhabited inductive predicates are not useful.

\begin{prop}[\textbf{EfQ}]
	\label{prop:efq}
	For any formula $A$, such that all inductive predicates occurring in $A$ are inhabited, we have that
	\[
	\big\{\bot\to X\,|\, X\in\mathbf{PP}(A)\big\}\vdash \bot\to A.
	\]
\end{prop}
In the following sections, we will often prove statements by \textit{simultaneous induction on formulas and predicates}. In such a proof, in the case of a predicate $\I\eqdef\pi_X\Iclauses{X}$ ($\pi=\mu,\nu)$, we will have to show that $\I$ has the desired property, assuming that all premises in all clauses i.e., $A_{im}(X)$ for $i<k,m<|\Vec{A}_i|$, already have this property. The only distinguishing feature of the formulas $A_{im}$, at this point, is $X\in\mathbf{SP}(A_{im})$, by which either $X\not\in$ $\mathbf{FP}(A_{im})$ i.e., $A_{im}$ is a parameter-premise, or $X\in\mathbf{PP}(A_{im})$ and all strictly-positive predicate-formulas $P\vec{t}$, where $X$ occurs freely, are of the form:
\[
\begin{aligned}
    &P\vec{t}=X\vec{t} &&\quad\mathrm{i.e.,\ }A_{im}\mathrm{\ is\ a\ recursive\ premise},\\
    &P\vec{t}=\vec{t}\in\J(X)&&\quad\mathrm{i.e.,\ }A_{im}\mathrm{is\ a\ nested\ premise\ with\ }\I\neq \J\in\mathcal{P}.
\end{aligned}
\]
This could be an unpleasant case-distinction, but we will avoid it in the next section with the \textit{monotonicity}-property, that we prove below. Its proof is straightforward, using induction on the definition of $\mathbf{SP}(\cdot)$, and it is omitted.

\begin{prop}[Monotonicity]
\label{prop:mon}
	Let $X\in\mathbf{SP}(A),\mathbf{SP(R})$ be strictly positive in a formula $A(X)$, respectively a predicate $R(X)$. If $P,Q\neq R$ are predicates, then
	\[
	\vdash (P\subseteq Q)\to (A[X/P]\subseteq A[X/Q]),\qquad\quad
        \vdash (P\subseteq Q)\to (R[X/P]\subseteq R[X/Q].
	\]
\end{prop}
\section{Strong Negation in TCF}
\label{sec: SNTCF}
We first motivate the definition of strong negation by giving a BHK-like explanation, that explains at the same time, how \textit{proofs}, respectively \textit{refutations}, can be combined to \textit{witness} the assertion or \textit{refutation} of compound statements.
\begin{description}[style=nextline,leftmargin=15pt,align=left,labelindent=5pt,]
	\item[\textbf{Conjunction}]
	A \textit{proof} of $A\wedge B$ is a \textit{proof} of $A$ together \textit{proof} of $B$.\\
	Any \textit{refutation} of either $A$ or $B$ \textit{refutes}  $A\wedge B$.
	\item[\textbf{Disjunction}]
	A \textit{proof} of $A\vee B$ is either a \textit{proof} of $A$ or $B$.\\
	A \textit{refutation} of $A$ together with a \textit{refutation} of $B$ refutes $A\vee B$.
	\item[\textbf{Universal}]
	A \textit{proof} of $\forall_x A$ is a routine assigning to any term $t$ a \textit{proof} of $A(t)$.
	A \textit{refutation} of $\forall_x A$ is a term $t$ together with a \textit{refutation} of $A(t)$.
	\item[\textbf{Existence}]
	A term $t$ with a \textit{proof} of $A(t)$ is a \textit{proof} of $\exists_x A$.\\
	A \textit{refutation} of $\exists_x A$ consists of a routine, which assigns a \textit{refutation} of $A(t)$ to any term \nolinebreak $t$.
	\item[\textbf{Implication}]
	A \textit{proof} of $A\to B$ is a routine, which assigns an proof of $A$ a proof of $B$.\\
	A \textit{proof} of $A$ together with a \textit{refutation} of $B$ \textit{refutes} $A\to B$.
\end{description}
\par
A logical system for constructive mathematics, which treats proof and refutations both as primitives, which are in some sense symmetric or dual, certainly sounds interresting. Here though, we prefer to stay inside the theory TCF i.e., a refutation is given by a proof of the strong negation. To clarify, we add one more rule namely,
\begin{description}[style=nextline,leftmargin=15pt,align=left,labelindent=5pt,]
	\item[\textbf{Strong Negation}]
	A \textit{proof} of the strong negation of $A$ is a \textit{refutation} of $A$.\\
	Any \textit{refutation} of $A$ is  \textit{refutes} a \textit{proof} of the strong negation of $A$
\end{description}
\par
To begin with, consider the following instructive example.
\begin{exa}[Accessible Part of a Relation]
	\label{ex:acc}
	If $\prec$ is a binary relation, then the accessible (from below) part of $\prec$ is inductively defined with the introduction rule
	\[
	\forall_x\big(\forall_y(y\prec x\to y\in \mathbf{Acc}_\prec)\to x\in\mathbf{Acc}_\prec\big)
	\]
	i.e., $x\in\mathbf{Acc}_\prec$ if any chain $x_0\prec\dots \prec x$ below $x$ is finite. Hence, $x\not\in\mathbf{Acc}_\prec$ means that not all chains below $x$ are finite. We can assert the failure of $x\in\mathbf{Acc}_\prec$ by exhibiting an infinite chain below $x$. Formally, $x\in \N{\mathbf{Acc}_\prec}$ exactly if there is $y\prec x$ with the same property. We introduce $\N{\mathbf{Acc}_\prec}$ as a coinductive predicate with closure rule
	\[
	x\in \N{\mathbf{Acc}_\prec}\to \exists_y\big(y\prec x\wedge y\in \N{\mathbf{Acc}_\prec}\big).
	\]
	The corresponding greatest fixed-point axiom is
	\[
		 (\N{\mathbf{Acc}_\prec})^\nu[Q]: \forall_x\big(x\in Q\to \exists_y(y\prec x\wedge y\in \N{\mathbf{Acc}_\prec}\cup Q)\big)\to Q\subseteq \N{\mathbf{Acc}_\prec}.
	\]
	It corresponds exactly to the inaccessible part of $\prec$, namely, $x\in \N{\mathbf{Acc}_\prec}$ if there is an infinite chain $\dots\prec y_1\prec y_0\prec x$. If for example, we consider the direct-predecessor order $\prec\eqdef \mu_X(\forall_n n\prec \mathtt{S}n)$ on natural numbers, then $n\in \mathbf{Acc}_\prec\leftrightarrow n\in\T{\Nat}$. Now we define the infinite natural number as the program constant
	\[
		\infty:\Nat,\qquad \infty=\mathtt{S}\infty.
	\]
	Particularly, $\infty\prec\mathtt{S}\infty=\infty$ and the premise of  $(\N{\mathbf{Acc}_\prec})^\nu[\{n\,|\,n\equiv\infty\}]$ is derivable i.e.,
	\[
		\vdash\exists_y(y\prec \infty\wedge y\in (\N{\mathbf{Acc}_\prec}\vee y\equiv\infty)).
	\]
	Hence,  $\infty\in\N{\mathbf{Acc}_\prec}$.
\end{exa}

Next, we define strong negation $A^{\No}$ of a formula $A$ of $\TCF$ by simultaneous recursion on
the definition of $\Form$ and $\Pred$.

\begin{defi}
	\label{defi:StrongNegation}
	To every formula $A(\vec{x})$ with free variables ranging over $\vec{x}$ we assign a formula $\N{A}(\vec{x})$ by the following rules:
	\[
	\begin{alignedat}{8}
		\N{(P\vec{t})}		&\eqdef\ \N{P}\vec{t},		&\quad
		&&\N{(A\vee B)}		&\eqdef\ \N{A}\wedge \N{B},&\\
		\N{(A\to B)}		&\eqdef\  A\wedge \N{B},	&\quad
		&&\N{(\forall_x A)}	&\eqdef\ \exists_x \N{A},	&\\
		\N{(A\wedge B)}		&\eqdef\ \N{A}\vee \N{B},	&\quad
		&&\N{(\exists_x A)}	&\eqdef\ \forall_x \N{A}.	&
	\end{alignedat}
	\]
	Assuming a fixed assignment $X\mapsto X',X'',X''',\dots$ of predicate variables to new distinct predicate variables of the same arity, we assign to each predicate $P$ of arity $\vec{\tau}$ a predicate $P^\mathbf{N}$ of the same arity by
	\[
	\N{X}\eqdef (X)', \qquad\N{\{\vec{x}\,|\, A\}}\eqdef\{\vec{x}\,|\,\N{A}\}
	\]
	and $(X^{(n)})'\eqdef X^{(n+1)}$. Furthermore, for two distinct predicate variables $X\neq Y$ we require $X^{(n)}\neq Y^{(m)}$
    For predicates $\I$ and $\J$ of the form
	\[
		\I\eqdef\mu_X\vec{K},\qquad \J\eqdef \nu_X\vec{K},\qquad \vec{K}(X)=\Iclauses{X},
	\]
	i.e., they are the least, respectively greatest, fixed-point of the same list of clauses,
	the definition of $\N{\I}$ and $\N{\J}$ are given by
	\[
	\begin{aligned}
		\N{\I}\eqdef\nu_{X'}\JNclause{X'},\\
		\N{\J}\eqdef\mu_{X'}\INclause{X'}.
	\end{aligned}
	\]
	Here, for $n_i:=\mathtt{lth}(\vec{A}_i)>0$ we have
	\[
		B_{i}:=\bigvee_{\nu<n_i}\N{A_{i\nu}}.
	\]
	If $n_i=0$ then $B_i \eqdef \mathtt{f}\equiv\mathtt{t}$ i.e., in this case we get the weak negation $\neg (\vec{y}\equiv\vec{t}_i)$.
\end{defi}
\begin{rem}
	\label{rem:sndef}
	Recall that the
	predicate $\I$ , as in the definition above, is equivalent to its \textit{operator-form}
	\[
		\vdash \I\vec{y}\leftrightarrow\bigvee_{i<k}\exists_{\vec{x}_i}\big(\vec{y}\equiv\vec{t}_i\wedge\bigwedge 	\vec{A}_i(\I)\big).
	\]
	Similarly, for $\N{\I}$, which has only one clause, we get
	\[
		\vdash\N{\I}\vec{y}\leftrightarrow\bigwedge_{i<k}\big(\forall_{\vec{x}_i}(\vec{y}\equiv\vec{t}_i\to \bigvee_{\nu<n_i}\N{A_{i\nu}}(\N{\I}))\big),
	\]
	and the right formula is \textit{almost} strong negation of the previous one on the right. In fact, they differ only by a classical implication in place of a material implication.
\end{rem}

\begin{rem}
    \label{rem:snindlog}
	Note that the logical connectives $\wedge,\vee,\exists$ do not need to be introduced as primitives and can be seen as special cases of inductive predicates. In the latter case, their strong negations according to Definition \ref{defi:StrongNegation} is consistent e.g.,
	disjunction and conjunction may be defined via
	\[
		A\vee' B\eqdef \mu_X(A\to X,B\to X),\qquad A\wedge'B\eqdef \mu_Y(A\to B\to Y).
	\]
	strong negations of these predicates are
	\[
		\N{(A\vee'B)}\eqdef \nu_{X'}(\N{A}\wedge'\N{B}\to X'),\qquad\N{(A\wedge'B)}\eqdef \nu_{Y'}(\N{A}\vee'\N{B}\to Y'),
	\]
	and since these are non-recursive we obtain
	\[
		\N{(A\vee'B)}\leftrightarrow \N{A}\wedge'\N{B},\qquad \N{(A\wedge'B)}\leftrightarrow\N{A}\vee'\N{B}.
	\]
\end{rem}

\begin{exa}[Even and odd]
	\label{ex:evenodd}
	We define even and odd numbers as the inductive predicates $\mathbf{Ev},\mathbf{Od}$, respectively, with the following introduction rules:
	\begin{gather*}
		\begin{aligned}
			&0\in\mathbf{Ev},\\
			&1\in\mathbf{Od},
		\end{aligned}\qquad
		\begin{aligned}
			&\forall_n(n\in\mathbf{Ev}\to (n+2)\in\mathbf{Ev}),\\
			&\forall_n(n\in\mathbf{Od}\to (n+2)\in\mathbf{Od}).
		\end{aligned}
	\end{gather*}
	Their strong negations are then coinductive predicates with the closure rules:
	\begin{gather*}
		\begin{aligned}
			&n\in\N{\mathbf{Ev}}\to n\not\equiv 0\wedge \forall_m\big(n\equiv(m+2)\to m\in\N{\mathbf{Ev}}\big),\\
			&n\in\N{\mathbf{Od}}\to n\not\equiv 1\wedge \forall_m\big(n\equiv(m+2)\to m\in\N{\mathbf{Od}}\big).
		\end{aligned}
	\end{gather*}
	With the least-fixed-point axioms of $\mathbf{Ev},\mathbf{Od}$ we easily verify that
	\[
	\forall_n\big(n\in\mathbf{Ev}\to n\in\N{\mathbf{Od}}\big),\qquad
	\forall_n\big(n\in\mathbf{Od}\to n\in\N{\mathbf{Ev}}\big).
	\]
	Furthermore, we  prove $\infty\in \N{\mathbf{Ev}},\N{\mathbf{Od}}$, since $\infty\neq \mathtt{0},\mathtt{1}$ and $\infty\equiv(\mathtt{SS}\infty)$.
	Hence, to prove the other direction we have to restrict to total natural numbers, and we get
	\[
	\forall_{n\in\T{\Nat}}\big(n\in\mathbf{Ev}\leftrightarrow n\in\N{\mathbf{Od}}\big),\qquad
	\forall_{n\in\T{\Nat}}\big(n\in\mathbf{Od}\leftrightarrow n\in\N{\mathbf{Ev}}\big).
	\]
\end{exa}

\begin{rem}
\label{rem:simneg}
    Although we can represent simultaneous least- and greatest fixed-points via nesting, the theory $\TCF$ may also be
    directly extended with \textit{simultaneous} definitions i.e.,
    \[
        (\I,\J)\eqdef \mu_{X,Y}\big[\Phi(X,Y)\subseteq X,\Psi(X,Y)\subseteq Y\big],
    \]
    where $X,Y$ must be strictly positive in both $\Phi,\Psi$ and with simultaneous least fixed-point axiom then is
    \[
        \Phi[\I\cap P,\J\cap Q]\subseteq P\to \Psi[\I\cap P,\J\cap Q]\subseteq P\to (\I,\J)\subseteq (P,Q).
    \]
    Definition \ref{defi:StrongNegation} can, more or less, be directly extended to cover simultaneous definitions. E.g., consider the simultaneous definition of Even and Odd
    \[
    \begin{aligned}
         (\mathbf{Ev}_{\mathrm{sim}},\mathbf{Od}_{\mathrm{sim}})\mueq\Big[&\big[0\in \mathbf{Ev}_{\mathrm{sim}},\forall_n(n\in \mathbf{Od}_{\mathrm{sim}} \to (n+1)\in \mathbf{Ev}_{\mathrm{sim}})\big],\\
         &\big[\forall_n(n\in \mathbf{Ev}_{\mathrm{sim}} \to (n+1)\in\mathbf{Od}_{\mathrm{sim}})\big]\Big].
    \end{aligned}
    \]
    Then strong negation will be
    \[
         \begin{aligned}
         (\N{\mathbf{Ev}}_{\mathrm{sim}},\N{\mathbf{Od}}_{\mathrm{sim}})\nueq
             \Big[&\big[\forall_n\big(n\not\equiv 0\wedge \forall_m(n\equiv (m+1)\to \N{m\in\mathbf{Od}}_{\mathrm{sim}})\to n\in\N{\mathbf{Ev}}_{\mathrm{sim}}\big)\big],
             \\
             &\big[\forall_n\big(\forall_m(n\equiv(m+1)\to m\in\N{\mathbf{Ev}}_{\mathrm{sim}})\to n\in\N{\mathbf{Od}}_{\mathrm{sim}}\big)\big]\Big].
         \end{aligned}
    \]
    Particularly, we  derive
    \[
        \begin{aligned}
        &\vdash\forall_k(k\equiv 1\to k\not\equiv 0\wedge\forall_m(k\equiv(m+1)\to m\in\N{\mathbf{Od}}_{\mathrm{sim}}\vee m\equiv 0)),\\
        &\vdash\forall_l(l\equiv 0\to \forall_m(l\equiv (m+1)\to m\in\N{\mathbf{Ev}}_{\mathrm{sim}}\vee m\equiv 0)),
        \end{aligned}
    \]
    since $k\equiv 1,k\equiv (m+1)\vdash m\equiv 0$ and $l\equiv 0,l\equiv (m+1)\vdash m\equiv 0$. For the latter we use transitivity, as $m\equiv \big[\mathrm{if}\ ((m+1)\neq 0)\ m\ 0\big]\equiv\big[\mathrm{if}\ (0\neq 0)\ m\ 0\big]\equiv 0$. Then we have derived the premise of the simultaneous greatest fixed-point axiom $(\N{\mathbf{Ev}}_{\mathrm{sim}},\N{\mathbf{Od}}_{\mathrm{sim}})^\nu\big[\{(k,l)\,|\, k\equiv 1\wedge l\equiv 0 \}\big]$, so that
    \[
        \vdash \big(1\in\N{\mathbf{Ev}}_{\mathrm{sim}}\big)\wedge 0\in\big(\N{\mathbf{Od}}_{\mathrm{sim}}\big).
    \]
\end{rem}

We prove, that strong negation of arithmetical falsity is derivable i.e., $\vdash \N{\bot}$. Furthermore, strong negation cancels weak negation.
\begin{lem}
	\label{lem:leibniz} We have that $\vdash\mathtt{f}\N{\equiv}\mathtt{t}$ and $\vdash\N{(\neg A)}\leftrightarrow A$.
\end{lem}
\begin{proof}
	Strong negation of \textit{Leibniz-equality} is
	\[
	\N{\equiv}\eqdef\mu_{X'}\big(\forall_{\vec{y}}(\forall_x(y_0\equiv x\wedge y_1\equiv x\to\bot)\to \vec{y}\in X')\big).
	\]
	As we can rename the bound predicate variables $X$ and $\equiv$ is non-recursive, we get
	\[
	\vdash (x\N{\equiv}y)\leftrightarrow \big(\mu_X\big(\forall_{\vec{y}}\forall_x(\vec{y}\not\equiv(x,x))\to \vec{y}\in X\big)\big)xy.
	\]
	By the transitivity of $\equiv$ we have $\mathtt{f}\equiv x\wedge \mathtt{t}\equiv x\to \mathtt{f}\equiv\mathtt{t}$, so $\forall_x \neg ((\mathtt{f},\mathtt{t})\equiv(x,x))$ and hence $\mathtt{f}\N{\equiv}\mathtt{t}$.
	Furthermore, $\vdash\N{(\neg A)}\leftrightarrow A$ since by definition $\N{(\neg A)}\eqdef A\wedge (\mathtt{f}\N{\equiv}\mathtt{t})$.
\end{proof}
\begin{rem}
	If we consider the equality of booleans $=_\Bol$ as the inductive predicate with introduction rules $\mathtt{f}\N{=_\Bol}\mathtt{f}$ and $\mathtt{t}\N{=_\Bol}\mathtt{t}\mueq (\mathtt{f}\N{=_\Bol}\mathtt{f},\mathtt{f}\N{=_\Bol}\mathtt{f})$, then its strong negation is
	\[
	\N{=_\Bol}\nueq \forall_{\vec{b}}\big((b_0\not\equiv\mathtt{f}\vee b_1\not\equiv\mathtt{f})\wedge (b_0\not\equiv\mathtt{t}\vee b_1\not\equiv\mathtt{t})\to b_0\N{=_\Bol}b_1\big).
	\]
	Thus we immediately get $(\mathtt{f},\mathtt{t}),(\mathtt{t},\mathtt{f})\in\N{=_\Bol}$. But this is not all, namely, if $b_0$ is provably a partial term i.e.,
	$\vdash b_0\not\equiv \mathtt{f}\wedge b_0\not\equiv \mathtt{t}$, then $(b,b_0),(b_0,b)\in \N{=_\Bol}$. Hence, $\N{=_\Bol}$ also contains all pairs of booleans where at least one of the components is partial. In the canonical model of $\TCF$ (see \cite{SW12}, Chapter 6) strong negation $\N{\mathbf{K}}$ of some inductive or coinductive predicate $\mathbf{K}$ can be seen as the weak complement of $\mathbf{K}\subseteq\{\vec{t}:\vec{\tau}\}$.
\end{rem}
\begin{exa}[Closure of a relation]
	\label{ex:closure}
	Let $\simeq$ be a binary relation. We define its reflexive, symmetric and transitive closure as the inductive predicate $\mathbf{Cl}(\simeq)$ with the introduction rules
	\begin{align*}
		&x\simeq y\to (x,y)\in \mathbf{Cl}(\simeq),\\
		&(x,x)\in \mathbf{Cl}(\simeq),\\
		&(x,y)\in \mathbf{Cl}(\simeq)\to (y,x)\in \mathbf{Cl}(\simeq),\\
		&(x,y)\in \mathbf{Cl}(\simeq)\to (y,z)\in \mathbf{Cl}(\simeq)\to (x,z)\in \mathbf{Cl}(\simeq).
	\end{align*}
	Its strong negation is then the greatest fixed-point with closure rule
	\begin{align*}
		(x,z)\in \N{(\mathbf{Cl}(\simeq))}\to \big(&(x\N{\simeq}z)\ \wedge\\
		&(x\not\equiv z)\ \wedge \\
		&((z,x)\in \N{(\mathbf{Cl}(\simeq))})\ \wedge \\
		&(\forall_y((x,y)\in \N{(\mathbf{Cl}(\simeq))}\vee (y,z)\in \N{(\mathbf{Cl}(\simeq))}))\big).
	\end{align*}
	From this closure rule we get that $\N{(\mathbf{Cl}(\simeq))}$ is an apartness relation. Moreover, we have that $\simeq\,\subseteq\mathbf{Cl}(\simeq)$ and $\N{(\mathbf{Cl}(\simeq))}\subseteq\, \N{\simeq}$ i.e.,
	\[
	\vdash x\simeq y\to(x,y)\in\mathbf{Cl}(\simeq),\qquad \vdash (x,y)\in\N{\mathbf{Cl}(\simeq)}\to x\N{\simeq}y.
	\]
	Hence the embedding $\simeq\subseteq \mathbf{Cl}(\simeq)$ is (provably) a strong implication in the sense of Rasiowa \cite{Ra74}.
	Even if $\simeq$ is already an equivalence relation particularly, $\mathbf{Cl}(\simeq)\subseteq\simeq$, then in general we cannot prove that $\N{\simeq}\subseteq\N{(\mathbf{Cl}(\simeq))}$. It remains to prove symmetry and co-transitivity  of $\N{\simeq}$. However, this does provide a method to define an apartness-relation $\#$ for any equivalence relation $\simeq\subseteq S\times S$, namely
	\[
	\#:=\N{(\mathbf{Cl}(\simeq))}_{S\times S}.
	\]
	Note that, by definition, $\#\subseteq \N{\simeq}$ however, in general we do not get the other direction. Particularly, $\#$ might be smaller than strong negation of $\simeq$.
\end{exa}

\section{Properties of Strong Negation in TCF}
\label{sec: propSNTCF}
Next we prove some fundamental properties of $A^{\No}$ in $\TCF$. In particular, we establish that $\N{A}$ behaves like a negation particularly,  that asserting $\N{A}$ is inconsistent with asserting $A$. Moreover, we prove that $\mathbf{NN}$-elimination holds.

The proofs will be given by induction on formulas, so we will also have to deal with the case of predicate variables. For the proofs to work, we will need to assume that the
predicate variable $X'$ associated to $X$ behaves like a strong (or weak) negation. We introduce the following notation:
\begin{defi}
	For $X$ of arity $\vec{\tau}$ and corresponding variables $\vec{x}$ we define the \textit{Ex-falso-} ($\mathbf{EF}$) and \textit{double-negation-elimination} ($\mathbf{DNE}$) properties

	\begin{gather*}
		\mathbf{EF}_{X;B}\eqdef\forall_{\vec{x}}(X\vec{x}\to X'\vec{x}\to B),\qquad \mathbf{DNE}_X\eqdef\forall_{\vec{x}}(X''\vec{x}\to X\vec{x}).
	\end{gather*}

	In the following we write $\lozenge_{\vec{X}}\vdash A$,
	where $\lozenge\in\{\mathbf{EF},\mathbf{DNE}\}$ and $\lozenge_{\vec{X}}\eqdef \{\lozenge_X\; |\; X\in\vec{X}\}$.
\end{defi}

\begin{rem}
\label{rem:dersub}
	In the following we will use substitutions for derivations several times. In general, if $\vec{A}\vdash B$ and $\eta$ is an \textit{admissible}, simultaneous substitution of type-,term-,predicate- and assumption-variables we have $\vec{A}\eta\vdash B\eta$. Particularly, if $\mathbf{EF}_{\vec{X};C}\vdash A$ and $B_0,B_1$ are formulas with
	$\vdash\forall_{\vec{x}}(B_0\to B_1\to C)$, equivalently $\vdash\mathbf{EF}_{X;C}[X\eqdef B_0,X'\eqdef B_1]$, then
	\[
		\vdash A[X\eqdef B_0,X'\eqdef B_1].
	\]
\end{rem}

\begin{lem}
	\label{lem:StrongNegationLem}
	For all formulas $A$ and $B$ we have that
	\[
	\mathbf{EF}_{\mathbf{PP}(A);B}\vdash \N{A}\to A\to B.
	\]
\end{lem}
\begin{proof}
	The proof is by  induction formulas and predicates.
	\begin{description}[style=unboxed,leftmargin=0mm,itemsep=5pt]
		\item[Case $A_0\to A_1$]
		We need to prove $$\mathbf{EF}_{\mathbf{PP}(A_0\to A_1);B}\vdash A_0\wedge \N{A_1}\to (A_0\to A_1)\to B.$$ Since $\mathbf{PP}(A_0\to A_1)=\mathbf{PP}(A_1)$, we apply the induction hypothesis for $A_1$.
		\item[Case $X$]
		 $X$ is a predicate variables and the goal by assumption $	\mathbf{EF}_{\{X\};B}\vdash \forall_{\vec{x}}(X\vec{x}\to X'\vec{x}\to B)$.
		\item[Case]
		For the remaining logical symbols we unfold definitions,
		\begin{alignat*}{3}
			\N{(A_0\wedge A_1)}&\eqdef \N{A}_0\vee \N{A}_1,&\qquad \N{(A_0\vee A_1)}&\eqdef \N{A}_0\wedge\N{A}_1,\\
			\N{(\forall_x A)}&\eqdef \exists_x \N{A},&\qquad \N{(\exists_xA)}&\eqdef \forall_x \N{A}.
		\end{alignat*}
		and can directly apply the induction hypothesis in all cases.
		\item[Case $\I\eqdef \mu_X(\Phi(X)\subseteq X)$] $\I$ is an Inductive predicate with
		\[
		\N{\I}\eqdef \nu_{X'}(X'\subseteq\N{\Phi}(X')),\qquad
		\mathbf{EF}_{\mathbf{PP}(\I);B}\vdash \Phi(X)\subseteq(\N{\Phi}(X')\to B).
		\]
		We instantiate $(\I)^\mu[P]$ with $P\vec{x}\eqdef \N{\I}\vec{x}\to B$, since the conclusion is exactly our goal
		\[
		\big(\Phi[\I\cap(\N{\I}\to B)]\subseteq (\N{\I}\to B)\big)\to \I\subseteq(\N{\I}\to B).
		\]
		In order to prove the premise, we unfold $\N{\I}=\N{\Phi}(\N{\I})$ in the premise directly above and it fits with the induction hypothesis. It only remains to derive the open assumption $\mathbf{EF}$, which is the tautology
		\[
		\I\cap(\N{\I}\to B)\to\N{\I}\to B.
		\]
		\item[Case] For a coinductive predicate $\J\eqdef\nu_Y(Y\subseteq\Psi(Y))$, its negation $\N{\J}\eqdef\mu_{Y'}(\N{\Psi}(Y')\subseteq Y')$ is a least fixed-point. We proceed exactly as in the last case, using $(\N{\J})^\mu$.\qedhere
	\end{description}
\end{proof}
\begin{rem}
	For the sake of clarity, we repeat the proof of Lemma \ref{lem:StrongNegationLem} for a specific predicate namely, total lists $\mathbf{L}(Y)\subseteq \mathbb{L}(\sigma)$ relative to $Y\subseteq\sigma$.
	This inductive predicate comes with axioms
	\[
	\begin{gathered}
		[\,]\in\mathbf{L}(Y),\qquad \forall_{y,l}\big(y\in Y\to l\in\mathbf{L}(Y)\to (y::l)\in\mathbf{L}(Y)\big),\\
		(\mathbf{L})^\mu[P]:[\,]\in P\to \forall_{y,l}\big(y\in Y\to l\in \mathbf{L}\to l\in P\to  (y::l)\in P\big)\to \mathbf{L}\subseteq P.
	\end{gathered}
	\]
	Moreover, the closure rule for $\N{\mathbf{L}}$ is given by
	\[
	l\in\N{\mathbf{L}}(Y')\to l\not\equiv[\,]\wedge\forall_{y,l_0}\big(l\equiv(y::l_0)\to y\in Y'\vee l_0\in \N{\mathbf{L}}(Y')\big).
	\] We aim to prove $\forall_x(Y'\to Y\to B)\vdash \forall_l(l\in\N{\mathbf{L}}\to l\in\mathbf{L}\to B)$, assume that $Y'\to Y\to B$ and use $(\mathbf{L})^\mu$ i.e., we need to prove $[\,]\in \N{\mathbf{L}}\to B$ and
	\begin{align*}
		\forall_{y,l}\big(y\in Y\to l\in\mathbf{L}\to (l\in \N{\mathbf{L}}\to B)\to(y::l)\in \N{\mathbf{L}}\to B\big).
	\end{align*}
	If $[\,]\in \N{\mathbf{L}}$, then $[\,]\not\equiv[\,]$ and we use \textbf{EfQ}.
	Otherwise, from $(y::l)\in \N{\mathbf{L}}$ we get $y\in Y'\vee l\in \N{\mathbf{L}}$.
	In the first case we use the assumption $Y'\to Y\to B$, and in the other cases we get $B$ by our assumptions.
\end{rem}
\begin{rem}
	Following Definition \ref{defi:StrongNegation}, double strong negation of an inductive or coinductive predicate is given by
	\begin{gather*}
		\NN{\I}\eqdef\mu_{X''}\big(\forall_{\vec{y}}(D\to X''\vec{y})\big),\qquad
		\NN{\J}\eqdef\nu_{X''}\big(\forall_{\vec{y}}(D\to X''\vec{y})\big),\\
		D\eqdef
		\bigvee_{i<k}\exists_{\vec{x}_i}\Big(\vec{y}\equiv\vec{t}_i\wedge \bigwedge_{\nu<n_i} \NN{A_{i\nu}}\Big).
	\end{gather*}
	This is equivalant to
	\[
	(\mu/\nu)_{X''}\big(\forall_{\vec{x}_i}((\NN{A_{i\nu}})_{\nu<n_i}\to X''\vec{t}_i)\big)_{i<k}.
	\]
\end{rem}
We now proceed to prove double negation elimination for strong negation.
\begin{thm}
	\label{thm:DoubleStrongNegElim}
	For all formulas $A$ we have that

	\[
	\mathbf{EF}_{\mathbf{NP}(A);\bot}\cup\mathbf{DNE}_{\mathbf{PP}(A)}\vdash\NN{A}\to A.
	\]
\end{thm}
\begin{proof}
	The proof is by induction on formulas and predicates.
	\begin{description}[style=unboxed,leftmargin=0mm,itemsep=5pt]
		\item[Case $A\to B$]
		 Then we have
		\[
		\mathbf{PP}(A)\subseteq\mathbf{NP}(A\to B),\qquad \mathbf{PP}(A\to B)=\mathbf{PP}(B).
		\]
		We have to prove
		\[
		\mathbf{EF}_{\mathbf{NP}(A\to B);\bot}\cup\mathbf{DNE}_{\mathbf{PP}(B)}\vdash\N{A}\vee\NN{B}\to A\to B.
		\]
		If $\N{A}$, then we use Lemma \ref{lem:StrongNegationLem} and \textbf{EfQ}. If $\NN{B}$, we apply the induction hypothesis.
		\item[Case $X$]
		$X$ is a predicate variable. We have to prove $\mathbf{DNE}_{X}\vdash \forall_{\vec{x}}(X''\vec{x}\to X\vec{x})$
		and assumption and goal coincide.
		\item[Cases]
		Other logical connectives unfold to
		\begin{alignat*}{3}
			\NN{(A\wedge B)}&\eqdef \NN{A}\wedge \NN{B},&\quad \NN{(A\vee B)}&\eqdef \NN{A}\vee\NN{B},\\
			\NN{(\forall_x A)}&\eqdef \forall_x \NN{A},&\quad \NN{(\exists_xA)}&\eqdef \exists_x \NN{A}.
		\end{alignat*}
		and the induction hypotheses are directly applicable.
		\item[Case $\J\eqdef\nu_Y(Y\subseteq\Psi(Y))$] $\J$ is a coinductive predicate. Then the definition of $\NN{\J}$ and the induction hypothesis are
		\[
		\NN{\J}\eqdef\mu_{Y''}(Y''\subseteq\NN{\Psi}(Y'')),\qquad \mathbf{DNE}_{Y}\vdash \NN{\Psi}(Y'')\subseteq\Psi(Y).
		\]
		For the goal we $\NN{\J}\subseteq\J$ we have to use $(\J)^\nu$, which presents the exercise $\NN{\J}\subseteq \Psi[\J\cup\NN{\J}]$. By unfolding $\NN{J}$ this follows from
		\[
		\NN{\Psi}[\NN{\J}]\subseteq\Psi[\J\cup\NN{\J}],
		\]
		a direct consequence of the induction hypothesis which, in turn requires proof of
		\[
		\mathbf{DNE}_{\J\cup\NN{\J}}:\NN{\J}\subseteq \J\cup\NN{\J}.
		\]
		In the case of an inductive predicate we proceed in a similar fashion.\qedhere
	\end{description}
\end{proof}
\begin{rem}
	For the sake of clarity, we redo the proof of Theorem \ref{thm:DoubleStrongNegElim} for the specific coinductive predicate representing nested trees with finite branching.
	The algebra $\mathbb{T}$ is generated by one constructor $\mathtt{Br}:\mathbb{L}(\mathbb{T})\to\mathbb{T}$ and we define the  (possibly) infinite, finitely branching trees as the coinductive predicate $\mathbf{Tn}\subseteq \mathbb{T}$ with the closure- respectively, greatest fixed-point rule
	\[
	\begin{aligned}
		(\mathbf{Tn})^-&: t\in \mathbf{Tn}\to \exists_l\big(l\in\mathbf{L}(\mathbf{Tn})\wedge t\equiv (\mathtt{Br}\ l)\big),\\
		(\mathbf{Tn})^\nu[P]&:\forall_t\big(Pt\to \exists_l(t\equiv(\mathtt{Br}\;l)\wedge l\in \mathbf{L}(\mathbf{Tn}\cup P))\big)\to P\subseteq \mathbf{Tn}.
	\end{aligned}
	\]
	The double strong negation is given then by
	\begin{align*}
		\NN{\mathbf{Tn}}\nueq \forall_{t}\big(t\in\NN{\mathbf{Tn}}\to \exists_{l}t\equiv (\mathtt{Br}\,l)\wedge l\in \NN{\mathbf{L}}[\NN{\mathbf{Tn}}]\big).
	\end{align*}
	We prove
	\[
	\vdash t\in \NN{\mathbf{Tn}} \to t\in\mathbf{Tn},
	\]
	and assume that we have already proved
	\[
	\mathbf{DNE}_Y\vdash l\in \NN{\mathbf{L}}(Y'')\to l\in \mathbf{L}(Y).
	\]
	By the greatest-fixed-point axiom of $\mathbf{Tn}$ it remains to show
	\[
	t\in\NN{\mathbf{Tn}}\to \exists_{l}\ \big(t\equiv\mathtt{Br}\,l\wedge l\in \mathbf{L}\big[\mathbf{Tn}\cup \NN{\mathbf{Tn}}\big]\big).
	\]
	By the closure-rule applied to $t\in\NN{\mathbf{Tn}}$ we obtain $l$ with $t\equiv (\mathtt{Br}\,l)$ and  $l\in \NN{\mathbf{L}}(\NN{\mathbf{Tn}})$. Since $\mathbf{DNE}_Y[Y:=\NN{\mathbf{Tn}},Y'':=\NN{\mathbf{Tn}}]$ holds trivially, we use the assumption to get $l\in \mathbf{L}(\NN{\mathbf{Tn}})$,
	and our goal $l\in \mathbf{L}(\mathbf{Tn}\cup \NN{\mathbf{Tn}})$ follows by monotonicity.
\end{rem}
\begin{rem}
We have established that the operation $\N{(\cdot)}$ shares properties with weak and classical negation and is \textit{strong} i.e., (disregarding predicate variables)
\begin{itemize}
	\item can express inconsistency $A,\N{A}\vdash \bot$ and hence is stronger than weak negation $\vdash \N{A}\to\neg A$,
	\item admits double-negation-elimination $\vdash\NN{A}\to A$,
	\item definitionally, $\N{(A\wedge B)}=\N{A}\vee\N{B}$ and $\N{\forall_x A}=\exists_x\N{A}$.
\end{itemize}
Next we will develop further properties of $\N{A}$. Particularly, we can easily exhibit counterexamples to the following schemas, usually associated with negations, either intuitionistic or classical
\[
	\begin{aligned}
		&\nvdash (\N{B}\to \N{A})\to A\to B,\quad&
		&\nvdash (A\to B)\to \N{B}\to \N{A},&\\
		&\nvdash A\to\NN{A},&
		&\nvdash \N{A}\to A^{\mathbf{NNN}},&\\
		&\nvdash (\N{A}\to A)\to A,&
		&\nvdash ((A\to B\N{)}\to C)\to A\to B\vee C.&
	\end{aligned}
\]
These stark differences to other negations are firstly due to the weakness of the implication-connective and the resulting, rather strong, $\N{(A\to B)}=A\wedge \N{B}$. Secondly, the \textit{involutive} part of the definition of strong negation, i.e., $\forall\slash\exists$ and $\wedge\slash\vee$, is also asymmetric, since $\exists\slash\vee$-statements tend to carry more information than $\forall\slash\wedge$-statements. Particularly, this explains, why proof by \textit{strong} contraposition is not admissible. Namely,
\[
	\nvdash (\N{B}\to \N{A})\to A\to B,
\]
is usually the case, if $B$ is a (non-trivial) $\exists$- or $\vee$-formula. Moreover, the other directions failure
\[
	\nvdash (A\to B)\to \N{B}\to \N{A},
\]
is most often due to $A$ being an $\to$-, $\forall$-, or $\wedge$-formula.
As an example take the real numbers $\mathbb{R}$ with equality $=_\mathbb{R}$ and apartness-relation $\#_\mathbb{R}=\N{(=_\mathbb{R})}$ and consider the following pair of strong contra-positives
\[
	\begin{aligned}
		&\vdash\forall_{x,y\in\mathbb{R}}(x\#_\mathbb{R}0\wedge y\#_\mathbb{R}0\to (x\cdot y)\#_\mathbb{R}0),\\
		&\nvdash\forall_{x,y\in\mathbb{R}}((x\cdot y)=_\mathbb{R}0\to x=_\mathbb{R}0\vee y=_\mathbb{R}0).
	\end{aligned}
\]
For a proof of the latter we need at least the \textit{Lesser Limited Principle of Omniscience} (\textbf{LLPO}) or, equivalently, the \textit{dichotomy} of reals
\[
	\forall_{x\in\mathbb{R}}\big(0\leq_{\mathbb{R}} x\vee x\leq_{\mathbb{R}}0\big),
\]
which is not derivable constructively.

In contrast to weak negation, which is nesting implications on the left, strong negation $\N{A}$ consumes implications. This is the reason for the failure $A\to\NN{A}$, which actually holds for implication-free formulas. Most of the negative statements above can be explained by the fact, that  proofs of $\N{A}$ underlie the same restriction as proofs of any other positive formula.
\end{rem}
As remarked above, if there are no implications, then $\N{A}$ simply flips $\wedge/\vee$ and $\forall/\exists$. Particularly $\NN{A}$ is an involution.
\begin{prop}
	\label{prop:noimpstab}
	For all implication-free formulas  $A$ (including predicates $\I$ and $\J$ with only implication-free premises) we have that
	\[
	\mathbf{Stab}_{\mathbf{PP}(A)}\vdash \forall_{\vec{y}}(A\to \NN{A}),
	\]
	where $\mathbf{Stab}_X\eqdef \forall_{\vec{x}}(X\vec{x}\to X''\vec{x})$.
\end{prop}
\begin{proof}
	By straightforward induction on formulas.
\end{proof}
\begin{lem}
	\label{thm:UnprovableEquiv} The following schemas are equivalent, where we assume that the formulas do not contain free predicate variables:
	\begin{description}[style=sameline,leftmargin=15pt,align=left,labelindent=5pt,]
		\item[$\mathbf{1}.$] $\vdash A\vee\N{A}$.
		\item[$\mathbf{2}.$] $\vdash (A\leftrightarrow B)\to (\N{A}\leftrightarrow \N{B})$.
		\item[$\mathbf{3}.$] $\vdash (A\to B)\to (\N{B}\to \N{A})$.
		\item[$\mathbf{4}.$]  $\vdash A\to \NN{A}$.
	\end{description}
\end{lem}
\begin{proof}
	We prove only the equivalence $\mathbf{1.}\Rightarrow\mathbf{4}$ and the other cases are equally as simple.
	\begin{description}[style=unboxed,leftmargin=0mm,itemsep=5pt]
		\item[Case $\mathbf{1.}\Rightarrow\mathbf{4}$]
		By assumption we have a derivation
		\[
		\vdash (A\to \NN{A})\vee (A\wedge A^{\mathbf{NNN}}).
		\]
		We proceed by case-distinction. In the second case  we have $A$ and $\N{A}$ by Theorem \ref{thm:DoubleStrongNegElim}. Then we can use Lemma \ref{lem:StrongNegationLem} to get $\NN{A}$.
		\item[Case $\mathbf{4.}\Rightarrow\mathbf{1}$]
		By assumption we have a derivation
		\[
		\vdash (A\to A)\to \N{A}\vee \NN{A}.
		\]
		Since $A\to A$ is derivable and $\NN{A}\to A$ by Theorem \ref{thm:DoubleStrongNegElim}, we get $\vdash A\vee \N{A}$.\qedhere
	\end{description}
\end{proof}
With this we can now show that there is no uniform way to define strong negation.
\begin{prop}
	\label{prop:nouniform}
	If there is a formula $P[X]$, such that for all formulas $A$
	\[
	\vdash \N{A}\leftrightarrow P[X\eqdef A],
	\]
	then $\vdash B\vee \N{B}$, for all formulas $B$ not containing free predicate variables.
\end{prop}
\begin{proof}
	Assume that a formula $P[X]$ as described above does exist. By induction on the structure of $P$, there is a derivation showing $$\vdash(A\leftrightarrow B)\to (P[X\eqdef A]\leftrightarrow P[X\eqdef B]),$$ for any two formulas $A,B$. But then, using the assumption $$\vdash (A\leftrightarrow B)\to (\N{A}\leftrightarrow \N{B}),$$ and Lemma \ref{thm:UnprovableEquiv} is applied.
\end{proof}
\begin{rem}
	As there are undecidable formulas, in general we have neither
	\[
	\vdash (A\to B)\to \N{B}\to\N{A}\qquad \text{nor}\qquad \vdash(A\leftrightarrow B)\to (\N{A}\leftrightarrow \N{B}).
	\]
	Constructively, the strong contrapositive $\N{B}\to\N{A}$ can be used to obtain more informative definitions. As an example consider the classical \textit{least-upper-bound} $lub(x,S)$ of a \textit{non-empty} and \textit{bounded} set $S\subseteq\mathbb{R}$
	\[
		lub(x,S)\eqdef \forall_y(\forall_z(z\in S\to z\leq y)\to x\leq y),
	\]
	provided $x$ is known to be an upper bound already.
	Namely, if $lub(x,S)$ then any other upper bound of $S$ is above $x$. Computationally this is not a very useful characterization, since it only tells us, that any upper bound is most definitely an upper bound. The \textit{strong contrapositive} however, is
	\[
		sup(x,S)\eqdef (\forall_y(y<x\to \exists_z(z\in S\wedge y<z)),
	\]
	the \textit{constructive} definition of the least upper bound and computationally much more useful.
\end{rem}
Hence, it is meaningful to investigate formulas $A,B$ that do have properties $2$ and $3$ from Lemma \ref{thm:UnprovableEquiv}. Following \cite{Ra74}, we incorporate strong implication and equivalence in TCF.
\begin{defi}
	\textit{Strong implication} $\overset{s}{\to}$ and strong equivalence $\overset{s}{\leftrightarrow}$ are defined by
	\[
	A\overset{s}{\to}B :\Leftrightarrow (A\to B)\wedge (B^\mathbf{N}\to A^\mathbf{N}),\qquad A\overset{s}{\leftrightarrow}B:\Leftrightarrow (A\leftrightarrow B)\wedge (A^\mathbf{N}\leftrightarrow B^\mathbf{N}).
	\]
\end{defi}
\begin{prop} If $\diamond\in\{\overset{s}{\to},\overset{s}{\leftrightarrow}\}$, then we get
\label{prop:strongimp}
		\[
		\begin{aligned}
		&\vdash(A\diamond B)\leftrightarrow(B^\mathbf{N}\diamond A^\mathbf{N}),\\
		&\vdash(B_0\diamond A_0)\to(A_1\diamond B_1)\to ((A_0\to A_1)\diamond (B_0\to B_1)).
		\end{aligned}
		\]
\end{prop}

\begin{proof}
	The proof is straightforward, we only prove one of the statements namely,

	\[
		(B_0\overset{s}{\to} A_0)\to(A_1\overset{s}{\to} B_1)\to ((A_0\to A_1)\overset{s}{\to} (B_0\to B_1)).
	\]
	It is our goal to derive
	\[
		(A_0\to A_1)\to B_0\to B_1,\qquad B_0\wedge\N{B_1}\to A_0\wedge \N{A_1},
	\]
	which, respectively, follow directly from the assumptions
	\[
		\left\{
			\begin{aligned}
				&B_0\to A_0,\\
				&\N{A_0}\to\N{B_0},
			\end{aligned}
		\right.\qquad
		\left\{
			\begin{aligned}
				&A_1\to B_1,\\
				&\N{B_1}\to\N{A_1}.
			\end{aligned}
		\right.\qedhere
	\]
\end{proof}

The \textit{stability} of a formula $A$ is denoted by $\mathbf{St}_A\eqdef \neg\neg A\to A$,
and the stability of a finite collections of formulae $\Gamma$ by $\mathbf{St}_\Gamma\eqdef\{\neg\neg B\to B\,|\, B\in\Gamma\}$. A formula $A$ is \textit{classically derivable}, $\vdash_{\mathrm{c}} A$, if and only if $\mathbf{St}_\Gamma\vdash A$ for a collection of formulae $\Gamma$. We know by Lemma \ref{lem:StrongNegationLem}, that $\vdash\N{A}\to\neg A$. We now prove that $\neg A$ and $\N{A}$ are equivalent with respect to classical derivability.
The proof is not complicated, but surprisingly, non-trivial. Due to the recursive definition of $\N{A}$, even for classical logic, an induction on formulae is needed to assert properties of $\N{A}$.

\begin{lem}
	\label{lem:classWNtoSN}
Let $A$ be a formula and $\vec{Y}=\mathbf{PP}(A)$, then
	\[\mathbf{Hyp}_{\vec{Y}}\vdash_{\mathrm{c}}\neg A\to\N{A},\]
where $\mathbf{Hyp}_{Y_i}\eqdef \neg X\to X'$.
\end{lem}
\begin{proof} The proof is by induction on formulas and predicates.
\begin{description}[style=unboxed,leftmargin=0mm,itemsep=5pt]
	\item[Case $X$] By assumption $\mathbf{Hyp}_X$.
	\item[Cases of symbols] We show
	\[
	\begin{gathered}
		\begin{alignedat}{3}
			\mathbf{Ih}_A&\vdash_{\mathrm{c}} \neg\exists_x \N{A}\to \forall_x A,&\qquad
			\mathbf{Ih}_{\{A,B\}}&\vdash_{\mathrm{c}} \neg(\N{A}\vee\N{B})\to A\wedge B,
			\\[3pt]
			\mathbf{Ih}_A&\vdash_{\mathrm{c}} \neg\exists_x A\to \forall_x \N{A},&
			\mathbf{Ih}_{\{A,B\}}&\vdash_{\mathrm{c}} \neg(A\vee B)\to \N{A}\wedge\N{B},
		\end{alignedat}\\[3pt]
		\mathbf{Ih}_{\{B\}}\vdash_{\mathrm{c}} \neg(A\to B)\to A\wedge \N{B},
	\end{gathered}
	\]
	where the open assumptions are
	$\mathbf{Ih}_A:\neg A\to\N{A}$ and $\mathbf{Ih}_B:\neg B\to\N{B}$.
	Then either the goal is reached directly or otherwise it follows with proof by contrapositive. We provide the following derivations, in natural deduction style.
	\[
	\begin{gathered}
	\begin{prooftree}
		\hypo{\ \mathbf{Ih}_A}
		\ellipsis{}{}
		\infer[no rule]1{\neg\N{A}\to A}
		\hypo{\ [\neg\exists_x\N{A}]}
		\hypo{[\N{A}]}
		\infer1{\ \exists_x\N{A}\ }
		\infer2{\bot}
		\infer1{\ \neg\N{A}\ }
		\infer2{A}
		\infer1{\ \neg\exists_x\N{A}\to \forall_x A\ }
	\end{prooftree}
	\qquad
	\begin{prooftree}
		\hypo{\ \mathbf{Ih}_A}
		\ellipsis{}{}
		\infer[no rule]1{\ \neg\N{A}\to A}
		\hypo{\ [\neg(\N{A}\vee\N{B})]}
		\hypo{[\N{A}]}
		\infer1{\ \N{A}\vee\N{B}\ }
		\infer2{\bot}
		\infer1{\ \neg \N{A}\ }
		\infer2{A}
		\hypo{\ \mathbf{Ih}_B}
		\ellipsis{}{}
		\infer[no rule]1{\ B\ }
		\infer2{A\wedge B}
		\infer1{\ \neg(\N{A}\vee\N{B})\to A\wedge B\ }
	\end{prooftree}
	\\[10pt]
	\begin{prooftree}
		\hypo{\ \mathbf{Ih}_A}
		\hypo{\ [\neg\exists_x A]}
		\hypo{[A]}
		\infer1{\ \exists_x A\ }
		\infer2{\bot}
		\infer1{\ \neg A\ }
		\infer2{\N{A}}
		\infer1{\ \neg\exists_x A\to \forall_x\N{A}\ }
	\end{prooftree}
	\qquad\qquad
	\begin{prooftree}[rule margin= 0.5ex]
		\hypo{\ \mathbf{Ih}_A}
		\hypo{\ [\neg(A\vee B)]}
		\hypo{[A]}
		\infer1{\ A\vee B\ }
		\infer2{\bot}
		\infer1{\ \neg A\ }
		\infer2{\ \N{A}}
		\hypo{\ \mathbf{Ih}_B}
		\ellipsis{}{}
		\infer[no rule]1{\ \N{B}\ }
		\infer2{\N{A}\wedge\N{B}}
		\infer1{\ \neg(A\vee B)\to \N{A}\wedge\N{B}\ }
	\end{prooftree}
\\[10pt]
	\begin{prooftree}
		\hypo{\ \mathbf{St}_A\hspace{-1em}}
		\hypo{\ [\neg(A\to B)]\hspace{-1.5em}}
		\hypo{\ \mathbf{EfQ}_{B}}
		\hypo{[A]}
		\hypo{[\neg A]}
		\infer2{\bot}
		\infer2{B}
		\infer1{\ A\to B\ }
		\infer2{\bot}
		\infer1{\ \neg\neg A\ }
		\infer2{A}
		\hypo{\ \mathbf{Ih}_B\hspace{-1em}}
		\hypo{\ [\neg(A\to B)]}
		\hypo{[B]}
		\infer1{\ A\to B\ }
		\infer2{\bot}
		\infer1{\ \neg B\ }
		\infer2{\N{B}\ }
		\infer2{A\wedge\N{B}}
		\infer1{\ \neg(A\to B)\to A\wedge\N{B}\ }
	\end{prooftree}
		\end{gathered}
	\]
	Note, that in the case of an implication the induction hypotheses is only used for the conclusion and $\mathbf{PP}(A\to B)=\mathbf{PP}(B)$.
	\item[Case $\I\eqdef \mu_X(\Phi(X)\subseteq X)$]
	$\I$ is inductive, so the induction hypothesis is
	\[
	\mathbf{Ih}_\I(X,X'):\neg X\to X'\vdash_{\mathrm{c}} \neg\Phi(X)\to \N{\Phi}(X').
	\]
	From this, we have to provide a derivation of $\neg \I\subseteq \N{\I}$,
	and since $\N{\I}\eqdef\nu_{X'}(X'\subseteq\N{\Phi}(X'))$, we may use the greatest fixed-point axiom $(\N{\I})^\nu[\neg\I]$ for that purpose. Then it suffices to show
	\[
	\neg\I\to\N{\Phi}[\N{\I}\cup\neg\I].
	\]
	We instantiate the induction hypothesis to
	\[
	\mathbf{Ih}_\I[\I,\N{\I}\cup\neg\I]\neg\I\to \N{\I}\cup\neg\I\vdash_{\mathrm{c}} \neg\Phi[\I]\to \N{\Phi}[\N{\I}\cup\neg\I].
	\]
	The open assumptions holds trivially and since $\I\leftrightarrow\Phi[\I]$, we also have $\neg\I\leftrightarrow \neg\Phi[\I]$ and hence
	\[
	\vdash_{\mathrm{c}} \neg\I\to \N{\Phi}[\N{\I}\cup(\neg\I)].
	\]
	\item[Case $\J\eqdef\nu_Z(Z\subseteq\Psi(Z))$] $\J$ is a coinductive predicate and we have the induction hypothesis
	\[
	\mathbf{Ih}_\J(Z,Z'):\neg Z\to Z'\vdash_{\mathrm{c}} \neg\Psi(Z)\to \N{\Psi}(Z').
	\]
	The goal is to find a derivation of $\neg \J\subseteq \N{\J}$. Since we presently find ourselves in classical logic, we may replace implications with their contrapositive and delete double weak negations on the way. Particularly, from the fixed-point axiom $(\J)^\nu[\neg\N{\J}]$ we get
	\[
	\vdash_{\mathrm{c}} \big(\neg\Psi[\J\cup\neg\N{\J}]\subseteq\N{\J}\big)\to\neg\J\subseteq\N{\J}.
	\]
	Moreover, by substitution in the induction hypothesis
	\[
	\mathbf{Ih}_\J[\J\cup\neg\N{\J},\N{\J}]:\neg \big(\J\cup\neg\N{\J}\big)\subseteq \N{\J}\vdash_{\mathrm{c}} \neg\Psi\big[\J\cup\neg\N{\J}\big]\to \N{\Psi}[\N{\J}].
	\]
	The assumption on the left is classically derivable via $\vdash\neg(\J\cup\neg\N{\J})\leftrightarrow\neg\J\cap\neg\neg\N{\J}$ and the stability $\mathbf{St}_{\N{\J}}$. Finally, $\N{\Psi}[\N{\J}]\subseteq\N{\J}$, and the proof is completed.\qedhere
\end{description}
\end{proof}

\section{Tight Formulas}

\label{sec: SNTCFtight}
In $\CM$ an inequality $\neq_X$ associated to a set $(X, =_X)$ is a binary relation, which is irreflexive and symmetric. It is called tight, if for every $x, x{'} \in X$ we have that
$\neg(x \neq_X x{'}) \to x =_X x{'}$. If the inequality $\neq_X$ is equivalent to the strong inequality on $X$
induced by $=_X$ i.e.,
$x \neq_X x{'} \TOT (x =_X x{'})^{\No}$, then $(X, =_X, \neq_X)$ is called a
\textit{strong} set, and then $=_X$ is called tight, if $\neq_X$ is tight. Similarly, a subset $A$ of a strong
set $X$ is called a \textit{tight subset}, if $\big(A^{\neq}\big)^c \subseteq A$. Within $\TCF$ we may now express the \textit{tightness} of a formula $A$ as
\[
	\vdash (\N{A}\to\bot)\to A.
\]

\begin{defi}
	We call a formula $A$ tight, if $\vdash (A^\mathbf{N}\to \bot)\to A$.
\end{defi}
In relation to strong implication, tight formulas are of interest, since they allow proof by strong contrapositive.
\begin{lem}
	Assume $A,B$ are formulae, then
	\[
		\mathbf{EF}_{\mathbf{PP}(A)}\cup\{\neg(\N{B})\to B\}\vdash (\N{B}\to\N{A})\to A\to B.
	\]
\end{lem}
\begin{proof}
	Assume $B$ is tight, $\N{B}\to\N{A}$ and $A$. By tightness of $B$, we further assume $\N{B}$ and prove $\bot$. The latter follows with $\mathbf{EfQ}$ and Lemma \ref{lem:StrongNegationLem} for the formula $A$.
\end{proof}
\begin{rem}
	As we are working in minimal logic, we are not especially attached to the arithmetical expression $\bot$. Due to the negative occurrence of $\bot$ in the definition of tightness above, we can replace it with something stronger and in fact, with any formula, inconsistent with $\N{A}$, particularly we have $\bot\vdash A$ and $\N{A},A\vdash \bot$, so
	\[
	\vdash(\N{A}\to\bot)\leftrightarrow(\N{A}\to A).
	\]
	Then $A$ is tight if and only if
	\[
		\vdash (\N{A}\to A)\to A.
	\]
\end{rem}
In the partial setting of \textbf{TCF} very few (non-trivial) formulas will be tight. The term \textit{tightness} is commonly used as an attribute of a relation and relations usually come with an intended domain. E.g., the pointwise equality on unary natural numbers is

\[
	\approx_\Nat\;\mueq\big[0\approx_\Nat 0,\forall_{n,m}(n\approx_\Nat m\to \mathtt{S}n\approx_\Nat\mathtt{S}m)\big],
\]
and the intended domain are the \textit{total} natural numbers $$\T{\Nat}\mueq[0\in\T{\Nat},\forall_n(n\in\T{\Nat}\to (\mathtt{S}n)\in\T{\Nat})].$$ The relation $\approx_\Nat$ is only going to be tight if we restrict appropriately namely,
\[
	\vdash\forall_{n,m\in\T{\Nat}}\big(\neg(n\N{\approx_\Nat}m)\to n\approx_\Nat m\big).
\]
In light of the discussion above, we introduce the notion of \textit{relative tightness}. Particularly, our notion will capture relations on cartesian products of finite, positive length i.e., we consider relativisations of the form
\[
	\forall_{x_0,\dots,x_{n-1}}\big(C_0(x_0)\to\dots\to C_{n-1}(x_{n-1})\to A\big),
\]
where $\mathbf{fv}(A)=\vec{x}$ are all the free variables of the formula $A$ and $\vec{C}$ are independent in the sense $\mathbf{fv}(C_i)=\{x_i\}$.
\begin{defi}
	\label{def:reltight}
	Let $A$ be a formula with free variables $\mathbf{fv}(A)=\vec{x}$ ($|\vec{x}| =n>0$) and $P$ a predicate of arity $\vec{\tau}$ ($|\vec{\tau}|=m>0)$. For $i<n$ let $C_i$ be a formula with $\mathbf{fv}(C_i)=\{x_i\}$. We call $A$ \textit{tight relative to} $\vec{C}$, $A\in\mathcal{T}_{\vec{C}}$, if
	\[
		\vdash \forall_{\vec{x}}\big(\vec{C}\to \neg \N{A}\to A\big).
	\]
    Furthermore, let $S_l$ be a predicate of arity $\tau_l$, for each $l<m$. We say that $P$ is \textit{tight relative} to $\vec{S}$, if
        \[
            \vdash \forall_{\vec{y}^{\vec{\tau}}}\big((y_l\in S_l)_{l<m}\to\neg(\vec{y}\in\N{P})\to \vec{y}\in P\big).
        \]
\end{defi}
Some basic properties are collected in the Proposition below.
\begin{prop}
	\label{prop:reltightbasic}
	Let $A,\vec{C}$ be formulae according to the requirements of Definition \ref{def:reltight} and, additionally assume $\mathbf{PP}(A)=\{\,\}$. Then the tightness of $A$ relative to $\vec{C}$ is located in between (strong) decidability and stability on $\vec{C}$, respectively i.e.,
	\[
		\begin{aligned}
			\vdash \forall_{\vec{x}}\big(\vec{C}\to A\vee\N{A}\big)\;\Rightarrow\; &A\in\mathcal{T}_{\vec{C}},\\
			&A\in\mathcal{T}_{\vec{C}} \;\Rightarrow\;\vdash\forall_{\vec{x}}\big(\vec{C}\to \neg\neg A\to A\big).
		\end{aligned}
	\]
	Moreover, if $A\in\mathcal{T}_{\vec{C}}$, then the formula $\forall_{\vec{x}}(\vec{C}\to A)$ is tight.
\end{prop}
\begin{proof}
	The first two are easy consequences of Proposition \ref{prop:efq} Lemma \ref{lem:StrongNegationLem}. For the third we need to prove
	\[
		\vdash \neg\N{\big(\forall_{\vec{x}}(\vec{C}\to A)\big)}\to \forall_{\vec{x}}(\vec{C}\to A).
	\]
	From the assumption
	\[
		\neg\N{\big[\forall_{\vec{x}}(\vec{C}\to A)\big]}\leftrightarrow \forall_{\vec{x}}\,\neg\Big[\Big(\!\bigwedge_{i<n}C_i\Big)\wedge \N{A}\Big],
	\]
	and $\vec{x}$ with $\vec{C}$, we obtain $\neg\N{A}$ and finally $A$ from $A\in\mathcal{T}_{\vec{C}}$.
\end{proof}
Next we characterize the interaction of the logical connectives with the notion of relative tightness. As anticipated, it behaves well in combination with $\to,\wedge,\forall$.
\begin{lem}
\label{lem:reltight}
	Let $A,B$ be formulas then
	\[
		\begin{gathered}
			A\in\mathcal{T}_{\vec{C}(\vec{x}),D(y)}\;\Rightarrow\; \big(\forall_{y} (D(y)\to A)\big)\in\mathcal{T}_{\vec{C}},\qquad\quad
			B\in\mathcal{T}_{\vec{D}(\vec{y})}
			\;\Rightarrow\;(A\to B)\in\mathcal{T}_{\vec{D}(\vec{y}),\vec{R}(\vec{x})},\\[5pt]
			\left.\begin{aligned}
				&A\in\mathcal{T}_{\vec{C}_A(\vec{x}),\vec{D}_A(\vec{y})},\\
				&B\in\mathcal{T}_{\vec{C}_B(\vec{x}),\vec{D}_B(\vec{z})}
			\end{aligned}\right\}\;\Rightarrow\;\big(A\wedge B\big)\in\mathcal{T}_{\vec{C}_A\wedge\vec{C}_B,\vec{D}_A,\vec{D}_B}
		\end{gathered}
	\]
	where in the second case, $\vec{x}\eqdef\mathbf{fv}(A)\setminus\mathbf{fv(B)}$ and for all $i<|\vec{x}|$ set $R_i\eqdef x_i\equiv x_i$  and, in the last case $\vec{C}_A\wedge\vec{C}_B=\big(C_{Ai}(x_i)\wedge C_{Bi}(x_i)\big)$.
\end{lem}
\begin{proof}
	\begin{description}[style=unboxed,leftmargin=0mm,itemsep=5pt]
		\item[Case 1] We need to prove
		\[
		\vdash \forall_{\vec{x}}\big(\vec{C}(\vec{x})\to \neg\big[\exists_y\big(D(y)\wedge \N{A}\big)\big]\to \forall_y(D(y)\to A)\big).
		\]
		We use the equivalence $\vdash\neg[\exists_y(D\wedge \N{A})] \leftrightarrow \forall_y \neg[D\wedge \N{A}],$
		which we can combine with the fact $\vdash \neg[D\wedge \N{A}]\to D\to \neg\N{A},$
		and the assumption.
		\item[Case 2]  	The goal is to find a derivation supporting
		\[
		\vdash \forall_{\vec{x},\vec{y}}\Big(\big[\vec{D}(\vec{y}),(x_i\equiv x_i)_{i<|\vec{x}|}\big]\to \neg\Big[A\wedge \N{B}\big]\to A\to B\big),
		\]
		from the assumption $B\in\mathcal{T}_{\vec{D}(\vec{y})}$ i.e.,
		\[
		    \vdash \forall_{\vec{y}}\big(\vec{D}(\vec{y})\to \neg\N{B}\to B\big).
		\]
		and follows in a similar way, using $\vdash \neg[A\wedge \N{B}]\to A\to \neg\N{B}$ and $(\equiv)^+$.
		\item[Case 3]
		The goal is unfolded to
		\[
		\forall_{\vec{x}}\big(\vec{E}\to \neg\big[\N{A}\vee\N{B}\big]\to A\wedge B\big]\big).
		\]
		Here we employ the fact $\vdash \neg[\N{A}\vee\N{B}]\leftrightarrow [(\neg\N{A})\wedge(\neg\N{B})]$.\qedhere
	\end{description}
\end{proof}

Lemma \ref{lem:reltight} can be seen as a recursive procedure generating tight formulas on a corresponding domain.
Hence, we now consider predicates $\I$ of arity $\vec{\tau}$ with $n\eqdef|\vec{\tau}|>0$,
\[
    \I\eqdef\omega_X\big[\forall_{\vec{x_i}}\big(\vec{A}_i(X)\to \vec{B}_i\to X\vec{t}_i\big)\big]_{i<k}\quad (\omega=\mu,\nu).
\]
In order to make an assumption $\neg[\vec{x}\in\N{\I}]$ useful, we introduce a reasonable domain-restriction for predicates.
\begin{defi}
    \label{def:tightproj}
    For $\I$ as above with $n>1$ and $l<n$ we say that $\I$ is \textit{non-overlapping}, if the terms $\vec{t}_i$ do not overlap i.e., for any non-constant map $f:\{0,\dots,n-1\}\to\{0,\dots,k-1\}$ and $i<k$ we have
\[
    \vdash \forall_{\vec{x_i},\vec{z}}\, \neg(\vec{t}_i\equiv(t_{(f\,m)m})_{m<k}),
\]
where $\vec{z}$ are all the free variables, which are not already contained in $\vec{x}_i$.
\end{defi}
\begin{rem}
Recall, that for $\I$'s  strong negation $\N{\I}$, we always have the property
\[
    \vdash \vec{s}\in\N\I\leftrightarrow \bigwedge_{i<k}\forall_{{x}_i}\big[\vec{s}\equiv\vec{t}_i\to\bigvee \N{\vec{A}_i\,\!\!}[X'/\N{\I}\big].
\]
Assume, that we have proven, that the terms $\vec{s}$ are not of the form of $\vec{t}_i$ ($i<k$) in the clauses of $\I$ i.e.,
\[
\vdash \bigwedge_{i<k}\forall_{\vec{x}_i}\vec{s}\not\equiv \vec{t}_i.
\]
Then we can prove $\vec{s}\in\N{\I}$, since all $i<k$ conjuncts follow with \textbf{EfQ}.\\
Usually it is not hard to prove relative tightness at all, but really, the goal is to find \textit{simple} and \textit{reasonable} $\vec{\phi}$ with $\I\in\mathcal{T}_{\vec{\phi}}$. Below, the last Lemma will provide relative tightness for very simple predicates akin to the pointwise equality $\approx_\Nat$ from the first example in the next section. Note that other domain-restriction, compared to non-overlapping, are possible and that the statement from below can be extended in various ways. Note however, that a formula, tight relative $\phi$ behaves almost classically in that domain. For many constructively interesting statements we can't find an acceptable restriction.
\end{rem}
\begin{lem}
    Let $\I$ be inductive, such that  all $i<k$ clauses are formed in the fashion
    \[
        \forall_{\vec{x_i}}\big((\vec{s}_{il}\in\I)_{l<n_i}\to \vec{t}_i\in\I\big),
    \]
    and furthermore, assume that $\mathbf{fv}(s_{ilm})\subseteq\mathbf{fv}(t_{im})$. Then, if $\I$ is also non-overlapping, we have $\I\in\mathcal{T}_{\vec{\J}}$, with
    \[
        \J_m\eqdef \mu_{Z_m}\big[\forall_{\vec{y_l}}(Z_m(s_{ilm})_{m<n_i}\to Z_m(t_{im}))\big]_{i<k}
    \]
\end{lem}
\begin{proof}
The predicates $\J_m$ ($m<n$) correspond exactly to the $m$-th coordinate-projections of $\I$. We have to show $y\in\I$ from $(y_m\in\J_m)_{m<n}$ and $\neg[\vec{y}\in\N{\I}$. We use all the axioms $(\J_m)^\mu$. Since $\neg[\vec{y}\in\N{\I}$, we can exclude, as in the remark above, all cases but $\vec{y}\equiv \vec{t}_{i}$. Then for each clause $l<n$ and $m<n_i$, the number of premises in $m$-th clause, we have a fitting induction hypotheses from $\J_m$
\[
    s_{ilm}\in \big(\J_m\cap\I\big)\quad (\mathrm{all\ m<n_i}).
\]
Particularly, we can apply $(\I)_i^+$.
\end{proof}
\section{Examples}
\label{sec: ex}
Next we consider some important case-studies and examples that reveal the naturality of Definition \ref{defi:StrongNegation}
and justify the use of $\TCF$ as a formal system for a large part of $\BISH$, as the use of strong negation
in $\TCF$ reflects accurately the seemingly ad hoc definitions of strong concepts in $\BISH$.

\begin{exa}[Similarity]
	\label{ex:sim}
	Assuming that we have some base type $\iota$, given by $i<k$ constructors $\mathtt{C}_i:\vec{\rho}_i\to\iota$, then equality for terms of type $\iota$ can be recursively characterized, namely $\mathtt{C}_i\vec{x}_i$ is equal to $\mathtt{C}_i\vec{y}_i$ whenever the lists $\vec{x}_i$ and $\vec{y}_i$ are.
	The equality of total terms of type $\iota$ is given by the similarity relation $\approx_\iota$, which is the inductive predicate with introduction rules
	\[
	\forall_{\vec{x},\vec{y}}\big((x_{\kappa}\approx y_\kappa)_{\kappa<\mathtt{lth}(\vec{\rho_i})}\to \mathtt{C}_i\vec{x}\approx_\iota\mathtt{C}_i\vec{y}\big),
	\]
	where $i<k$, and $(x_{\kappa}\approx y_\kappa)$ uses the corresponding similarity for the type $\rho_{i\kappa}$.
	Then, $\N{(\approx_{\iota})}$ is the greatest fixed-point with closure-rule
	\begin{align*}
		t\N{(\approx_{\iota})}s\to
		\bigwedge_{i<k}\bigg(\forall_{\vec{x}_i,\vec{y}_i}\Big((t,s)\equiv(\mathtt{C}_i\vec{x}_i,\mathtt{C}_i\vec{y}_i)\to\bigvee_{\kappa}\N{(x_{i\kappa}\approx y_{i\kappa})}\Big)\bigg).
	\end{align*}
	For the natural numbers the introduction, respectively closure, rules are given by
	\begin{gather*}
		0\approx_\Nat 0,\qquad \forall_{n,m}(n\approx_\Nat m\to \mathtt{S}n\approx_\Nat\mathtt{S}m),\\
		n \N{\approx_{\Nat}} m \to (n,m)\not\equiv(0,0)\wedge
		\forall_{k,l}\big((n,m)\equiv (\mathtt{S}k,\mathtt{S}l)\to k\N{\approx_{\Nat}}l\big).
	\end{gather*}
	If we restrict to $\T{\Nat}$, this is a tight apartness relation since then $\approx_\Nat$ is decidable. Moreover, any total function $f:\Nat\to\Nat$ is extensional and strongly extensional
	\[
	\forall_{n,m\in\T{\Nat}}\big(n\approx_{\Nat}m\to fn\approx_{\Nat}fm\big),\qquad \forall_{n,m\in\T{\Nat}}\big(fn\N{\approx_\Nat}fm\to n\N{\approx_\Nat}m\big),
	\]
	where the latter follows from the derivation
	\[
	\vdash \forall_{n,k\in\T{\Nat}}(fn\N{\approx}f(n+k)\to n\N{\approx_\Nat}(n+k)).
	\]
	Although the case of natural numbers is easy, for any other finitary algebra $\iota$ we can proceed in a similar manner. If we consider types of higher level e.g., functionals $F$ of type $(\Nat\to\Nat)\to\Nat$, then in general we can't derive the following implication
	\begin{align*}
		\forall_{f}&\big(\forall_{n\in\T{\Nat}}(fn\in\T{\Nat})\to Ff\in\T{\Nat}\big)\to\\ &\forall_{f,g}\big(\forall_{n,m}(n\approx m\to fn\approx gn)\to Ff\approx Fg\big).
	\end{align*}
	However, for any particular $F$ with a given proof term $M$ of  $\forall_{f}\big(\forall_{n\in\T{\Nat}}(fn\in\T{\Nat})\to Ff\in\T{\Nat}\big)$, its extracted term $et(M)$ (see \cite{SW12}, Section 7.2.5) is continuous, which is enough to prove strong extensionality of $F$.
\end{exa}

\begin{exa}[Bisimilarity and Continuity]
	\label{ex:bisim}
	Instead of working with an equality on the total terms of a base-type $\iota$, we might consider an equality on the cototal (possibly infinite) terms. If we take the clauses of similarity $\approx_\iota$ from the last example and consider the greatest instead of the least fixed-point we get the bisimilarity relation $\sim_\iota$.
	For streams of booleans $\mathbb{S}$ its bisimilarity $\sim_{\mathbb{S}}$ is defined as the greatest fixed point with closure rule
	\[
	u_0\sim_{\mathbb{S}} u_1 \to \exists_{\vec{b},\vec{v}}(b_0=_\Bol b_1 \wedge v_0\sim_{\mathbb{S}} v_1).
	\]
	Its strong negation is the least-fixed-point with introduction rule
	\[
	\forall_{\vec{b},\vec{v}}\big(\vec{u}\equiv \vec{b}::\vec{v}\to b_0\N{=_\Bol}b_1\vee v_0\N{\sim_{\mathbb{S}}}v_1\big)\to u_0\N{\sim_{\mathbb{S}}}u_1.
	\]
	i.e., $u_0\N{\sim_{\mathbb{S}}}u_1$ exactly if after some finite time $n\in\Nat$ of comparing $u_0$ and $u_1$ we find two digits $(u_0)_n\N{=_\Bol}(u_1)_n$. Let $h$ be a cototal function of type $\mathbb{S}\to \mathbb{B}$ i.e., $\forall_u(u\sim_{\mathbb{S}} u\to hu\in\T{\Bol})$. If $h$ is equipped with a proof-term of its cototality we can prove its strong extensionality, but in general we have to assume that $h$ is also continuous i.e.,
	\[
	\forall_u \exists_n \forall_v h(u)=_\Bol h(u|_n\ast v),
	\]
	where $u|_n$ is the initial segment of $u$ of length $n$ and $\ast$ is the concatenation of a list of booleans to a stream. From the continuity hypothesis we get the extensionality of $h$, since $u\sim_{\mathbb{S}} v$ implies that all initial segments $u|_n$ and $v|_n$ are equal lists. By taking the maximum of the continuity moduli of $u$ and $v$, we conclude that $h(u)=_\Bol h(v)$. Moreover, we also get strong extensionality of $h$ namely,
	\[
	\forall_{u,v}\big(hu\N{=_\Bol}hv\to u\N{\sim_{\mathbb{S}}}v\big).
	\]
	This is the case, because we can prove
	\[
	\forall_{u,v,n,w}\big(h(u|_n\ast w)\N{=_\Bol}h(v|_n\ast w)\to (u|_n\ast w)\N{\sim_{\mathbb{S}}}(v|_n\ast w)\big)
	\]
	by induction on $n$ as follows. If $n=0$, then $h(w)\N{=_\Bol}h(w)$ contradicts the irreflexivity of $\N{=_\Bol}$. In the step case either $(u)_{n+1}=_\Bol (v)_{n+1}$, and we use the induction hypothesis, or $(u)_{n+1}\N{=_\Bol}(v)_{n+1}$. Moreover, if $h$ is non-extensional, then it is discontinuous i.e.,
	\[
	\exists_{u,v}\big(u\sim_{\mathbb{S}}v\wedge hu\N{=}hv\big)\to \exists_u \forall_n \exists_v h(u|_n\ast v)\N{=_\Bol} hu.
	\]
	Namely, if $h$ is non-extensional for $u,v$ then $h(u|_n\ast {_n{|}}v)\N{=_\Bol}hu$, where $v\equiv v|_n\ast {_n{|}}v \equiv u|_n\ast {_n{|}}v$.
\end{exa}

\begin{exa}[Cauchy reals, equality and apartness]
	\label{ex:real}
	We consider the Cauchy-reals of type $\Rea:=(\Nat\to\mathbb{Q})\times(\mathbb{\mathbb{P}\to\Nat})$ i.e., rational sequences with an explicit modulus of Cauchyness, where $\mathbb{P}$ is the type of positive integers. Furthermore, we restrict to $x,y$ of type $\Rea$ that are total, namely, $x\in\T{\Rea}\Leftrightarrow x=(x_s,x_m)$, where $x_s,x_m$ are total functions, and consider pairs $(x_s,x_m)$ that are Cauchy-sequences i.e.,
	\[
	\forall_{p}\forall_{n,k\geq (x_mp)}\big(|(x_s\,n)-(x_s\,k)|\leq_\mathbb{Q} 2^{-p}\big).
	\]
	Nonnegativity and positivity are defined by
	\begin{gather*}
		0\leq_\Rea x \eqdef \forall_{p}\exists_{n\geq (x_mp)}\big(-2^{-p}<_{\mathbb{Q}} (x_s\,n)\big),\\ 0<_\Rea x\eqdef \exists_{p}\forall_{n\geq (x_m\,p)} \big(2^{-p}\leq_{\mathbb{Q}} (x_s\,n)\big),
	\end{gather*}
	respectively. The less ($<_\mathbb{Q}$) and less-than-or-equal ($\leq_\mathbb{Q}$) relation for rational numbers are terms of type $\mathbb{Q}\to\mathbb{Q}\to\mathbb{B}$, which are decidable on total inputs. In particular, $q_0\leq_{\mathbb{Q}} q_1=\mathtt{b}\leftrightarrow q_1<_{\mathbb{Q}}q_0=\mathtt{\overline{b}}$. Hence strong negation of nonnegativity is
	\[
	\N{(0\leq_\Rea x)}\eqdef \exists_{p}\forall_{n\geq (x_m p)}\big(-2^{-p}>_{\mathbb{Q}} (x_s\,n)\big).
	\]
	Then we immediately get $\N{(0\leq_\Rea x)}\leftrightarrow (0<_\Rea-x)$. Notice that due to the Cauchyness of $(x_s,x_m)$ we get
	\begin{gather*}
		\vdash (0\leq_\Rea x)\leftrightarrow \forall_p \big(-2^{-p}\leq (x_s(x_m\, p))\big),\\ \vdash (0<_\Rea x)\leftrightarrow \exists_p \big(2^{-p}\leq (x_s(x_m(p+1)))\big).
	\end{gather*}
	If we define
	\begin{gather*}
		x\leq_\Rea y\eqdef 0\leq_\Rea y-x,\qquad x<_\Rea y\eqdef 0<_\Rea y-x,\\
		x=_\Rea y\eqdef x\leq_\Rea y\wedge y\leq_\Rea x,
	\end{gather*}
	then strong negation of equality on reals is
	\[
	\N{(x=_\Rea y)}\eqdef y<_\Rea x\vee x<_\Rea y,
	\]
	which is the canonical tight apartness relation on $\Rea$. We also write $x<_p y$, if we explicitly need the witness $p$. If we define the absolute value $|(z_s,z_m)|\eqdef ((|z_s\,n|)_n,z_m)$, then we get
	\[
	x\N{=_\Rea} y\leftrightarrow \exists_{p\in\T{\mathbb{P}}}( 0<_p |x-y|).
	\]
\end{exa}

\begin{exa}[Real functions]
	\label{ex:realfun}
	Now we consider functions $f$ of type $\Rea\to\Rea$. We assume that every $f$ has a codomain $I\subseteq \Rea$, such that $f$ is total on $I$ and maps Cauchy-sequences in $I$ to Cauchy-sequences in $\mathbb{R}$.
	Using the definitions in Example \ref{ex:real} we consider monotone, injective, extensional functions and their respective strong versions:
	\begin{gather*}
		\forall_{x,y\in I}(x\leq_\Rea y\to fx\leq_\Rea fy),\qquad \forall_{x,y\in I}(fy<_\Rea fx\to y<_\Rea x),\\
		\forall_{x,y\in I}(fx=_\Rea fy\to x=_\Rea y),\qquad \forall_{x,y\in I}(x\not=_\Rea y\to fx\not=_\Rea fy),\\
		\forall_{x,y\in I}(x=_\Rea y\to fx=_\Rea fy),\qquad \forall_{x,y\in I}(fx\not=_\Rea fy\to x\not=_\Rea y).
	\end{gather*}
	The existence of a derivation showing
	\[
	\vdash\forall_{x,y\in I}(x\leq_\Rea y\to fx\leq_\Rea fy)\to\forall_{x,y\in I}(fy<_\Rea fx\to y<_\Rea x),
	\]
	is equivalent to Markov's principle. To show this let $x,y$ with $\neg(y\leq_\Rea x)$. and define $f$ with $I\eqdef\{x,y\}$ by
	$f(x)\eqdef 0$ and $f(y)\eqdef 1$. Then $f$ is monotone and $fx<_\Rea fy$, hence $x<_\Rea y$. The same can be proven for strong injectivity and extensionality, respectively. A function $f$ is (uniformly) continuous, if there exists a total and monotone modulus $M$ of type $\mathbb{P}\to\Nat$ with
	\[
	\forall_{x,y\in I,p}\big(|x-y|\leq_\Rea 2^{-(Mp)}\to |fx-fy|\leq_\Rea 2^{-p}\big).
	\]
	Then we can show that continuous functions $f$ are (strongly) extensional and the monotonicity and the injectivity hypothesis imply their respective strong versions. It is easy to prove extensionality from continuity. For strong extensionality we need to show
	\[
	\forall_{x,y}(0<_\Rea|fx-fy|\to 0<_\Rea |x-y|).
	\]
	For simplicity, assume $y=f(y)=0$. Then by the assumption $0<_\Rea |fx|$ there is $p_0$ with $\forall_{n\geq(|fx|_mp_0)} (2^{-p_0}\leq_{\mathbb{Q}} (|fx|_sn))$. We have to show $\exists_p\forall_{n\geq (x_mp)} (2^{-p}\leq_{\mathbb{Q}} (|x|_s\,n))$, for which we take $p\eqdef Mp_0+1$. Assume $n\geq x_m(M(p_0)+1)$ and $|x|_sn<_{\mathbb{Q}} 2^{-(Mp_0+1)}$. Then $|x|\leq_\Rea 2^{-Mp_0}$ and by the continuity assumption we get $|fx|\leq_\Rea 2^{-p_0}$. By definition this is
	\[
	\forall_p\exists_{n\geq (|fx|_mp)} ((|fx|_s\,n)\leq_{\mathbb{Q}}2^{-p_0}+2^{-p}),
	\]
	contradicting the initial assumption for $p\eqdef p_0$. Hence, $2^{-(Mp_0+1)}\leq_{\mathbb{Q}} |x|_sn$ and $0<_\Rea |x|$.
\end{exa}

\section{Concluding comments}
\label{sec: concl}

In this paper we incorporated strong negation in the theory of computable functionals $\TCF$ by defining
strong negation $A^{\No}$ of a formula $A$ and strong negation $P^{\No}$ of a predicate $P$ in $\TCF$.
As far as we know,
the extension of strong negation to inductive and coinductively defined predicates is new.
The fact that strong negation of an inductively defined predicate is an appropriate coinductively defined one,
and vice versa, reveals a new aspect of the duality between induction and coinduction.
We also presented some basic properties
of strong negation in $\TCF$ and, in analogy to $\CM$, we introduced the so-called tight formulas and relative tight formulas of $\TCF$.

Examples \ref{ex:evenodd} and \ref{ex:acc} reflect the ``naturality'' of our main Definition \ref{defi:StrongNegation}. In Example \ref{ex:evenodd} the inductive predicates $\mathbf{Ev}$
and $\mathbf{Od}$ are naturally connected through their strong negation, while in Example \ref{ex:acc}
strong negation of the inductive predicate ``accessible part of a relation'' corresponds to its inaccessible part.

Examples \ref{ex:closure}--\ref{ex:realfun} reinforce the appropriateness of $\TCF$ as a formal system for Bishop-style constructive mathematics.
According to Example \ref{ex:closure}, strong negation of an equality is an apartness relation, while in Example \ref{ex:sim}
strong negation of the similarity relation of total natural numbers is shown to be a tight apartness relation
and the total functions of type $\Nat \to \Nat$ are both extensional and strongly extensional with respect to the similarity relation of
naturals. The relation of bisimilarity to continuity is explained in Example \ref{ex:bisim}. In Example \ref{ex:real} strong negation of the equality of reals is shown to be
the canonical apartness relation on reals. In Example \ref{ex:realfun} known facts about strongly extensional functions of type $\Rea\to\Rea$ are recovered within TCF.

Moreover, all introduced notions, examples, and results can be easily formalised in the proof assistant Minlog~\cite{TMS} that accommodates TCF.

According to Feferman~\cite{Fe79}, a formal theory $T$ is \textit{adequate} for an informal body of mathematics $M$,
if every concept, argument, and result of $M$ is represented by a (basic or defined) concept, proof,
and a Theorem, respectively, of $T$. Additionally, $T$ is \textit{faithful} to $M$, if
every basic concept of $T$ corresponds to a basic concept of $M$ and every
axiom and rule of $T$ corresponds to or is implicit in the assumptions and reasoning followed in $M$ i.e., $T$
does not go beyond $M$ conceptually or in principle. The use of strong negation in $\TCF$, and especially in its ``total''
fragment, provides extra evidence for its adequacy as a formal system for $\BISH$. As dependent operations are crucial to
$\BISH$ (see~\cite{Pe20}), one needs to add dependency to $\TCF$, in order to get an adequate formalisation of $\BISH$.
Such an extension of $\TCF$ is an important open problem of independent interest. The use of inductive and
coinductive predicates in $\TCF$ prevents $\TCF$ from being a faithful formalisation\footnote{Equality is also treated differently in $\TCF$. Moreover, partiality in $\BISH$ is defined through totality, while in $\TCF$ partiality is a primitive concept.}\footnote{The presence of impredicativity in~\cite{Bi67}, and even more in~\cite{BB85}, is addressed in~\cite{Pe20} and~\cite{Pe22a}, where set-indexing methods and witnessing data are added, in order to bypass the aforementioned impredicativity. Moreover, TCF is predicative and this fact justifies the use thereof for formalising BISH.} for $\BISH$.
Although inductive definitions with rules of countably  many premises are included in the extension\footnote{As we have remarked already in the introduction, $\BISH^*$
	is the system corresponding to~\cite{Bi67} and to the constructive topology of Bishop spaces (see~\cite{Pe15, Pe19b, Pe20b}
	and~\cite{Pe23}), while $\BISH$ is the system corresponding to~\cite{BB85}.}
$\BISH^*$ of $\BISH$, as far as we know there is no single case-study of a coinductive definition within BISH or BISH$^*$.
\bibliography{library}
\bibliographystyle{alphaurl}
\end{document}